\documentclass[final,onefignum,onetabnum]{siamart220329}

\usepackage{braket,amsfonts}
\usepackage{amsmath,empheq}
\usepackage{cleveref}
\usepackage{mathtools}
\usepackage{array}

\usepackage{caption}
\usepackage{subcaption}

\usepackage{pgfplots}

\newsiamthm{claim}{Claim}
\newsiamremark{remark}{Remark}
\newsiamremark{hypothesis}{Hypothesis}
\crefname{hypothesis}{Hypothesis}{Hypotheses}

\usepackage{algorithmic}

\usepackage{graphicx,epstopdf}

\Crefname{ALC@unique}{Line}{Lines}

\usepackage{amsopn}

\usepackage{xspace}
\usepackage{bold-extra}
\usepackage[most]{tcolorbox}

\colorlet{texcscolor}{blue!50!black}
\colorlet{texemcolor}{red!70!black}
\colorlet{texpreamble}{red!70!black}
\colorlet{codebackground}{black!25!white!25}

\newcommand{\E}{\mathcal{E}}
\newcommand{\norm}[1]{\left\|#1\right\|}
\newcommand{\inner}[1]{\left\langle #1\right\rangle}
\newcommand{\brak}[1]{\left( #1\right)}

\lstdefinestyle{siamlatex}{%
  style=tcblatex,
  texcsstyle=*\color{texcscolor},
  texcsstyle=[2]\color{texemcolor},
  keywordstyle=[2]\color{texemcolor},
  moretexcs={cref,Cref,maketitle,mathcal,text,headers,email,url},
}

\tcbset{%
  colframe=black!75!white!75,
  coltitle=white,
  colback=codebackground, % bottom/left side
  colbacklower=white, % top/right side
  fonttitle=\bfseries,
  arc=0pt,outer arc=0pt,
  top=1pt,bottom=1pt,left=1mm,right=1mm,middle=1mm,boxsep=1mm,
  leftrule=0.3mm,rightrule=0.3mm,toprule=0.3mm,bottomrule=0.3mm,
  listing options={style=siamlatex}
}

\newtcblisting[use counter=example]{example}[2][]{%
  title={Example~\thetcbcounter: #2},#1}

\newtcbinputlisting[use counter=example]{\examplefile}[3][]{%
  title={Example~\thetcbcounter: #2},listing file={#3},#1}

\DeclareTotalTCBox{\code}{ v O{} }
{ %fontupper=\ttfamily\color{texemcolor},
  fontupper=\ttfamily\color{black},
  nobeforeafter,
  tcbox raise base,
  colback=codebackground,colframe=white,
  top=0pt,bottom=0pt,left=0mm,right=0mm,
  leftrule=0pt,rightrule=0pt,toprule=0mm,bottomrule=0mm,
  boxsep=0.5mm,
  #2}{#1}

\patchcmd\newpage{\vfil}{}{}{}
\begin{tcbverbatimwrite}{tmp_\jobname_header.tex}
\title{Linear convergence of forward-backward accelerated algorithms without knowledge of the modulus of strong convexity \thanks{\funding{This work was supported by Grant No.YSBR-034 of CAS and Grant No.12288201 of NSFC.}}}

\author{Bowen Li\footnotemark[2] \and Bin Shi\thanks{Academy of Mathematics and Systems Science, Chinese Academy of Sciences, Beijing 100190, China (\email{shibin@lsec.cc.ac.cn}); \; School of Mathematical Sciences, University of Chinese Academy of Sciences, Beijing 100049, China} \and Ya-xiang Yuan\footnotemark[2]}

\headers{Linear convergence of forward-backward accelerated algorithms without knowledge of the modulus of strong convexity}{ Bowen Li and Bin Shi and Ya-xiang Yuan}
\end{tcbverbatimwrite}
\title{Linear convergence of forward-backward accelerated algorithms without strongly convex modules\thanks{\funding{This work was supported by Grant No.YSBR-034 of CAS and Grant No.12288201 of NSFC.}}}

\author{Bowen Li\footnotemark[2] \and Bin Shi\thanks{Academy of Mathematics and Systems Science, Chinese Academy of Sciences, Beijing 100190, China (\email{shibin@lsec.cc.ac.cn}); \; School of Mathematical Sciences, University of Chinese Academy of Sciences, Beijing 100049, China} \and Ya-xiang Yuan\footnotemark[2]}

\headers{Linear convergence of forward-backward accelerated algorithms without strongly convex modules}{ Bowen Li and Bin Shi and Ya-xiang Yuan}

\ifpdf
\hypersetup{ pdftitle={Linear convergence of forward-backward accelerated algorithms without strongly convex modules} }
\fi

\begin{document}
\maketitle

%% ------------------------------------------------------------------
%% ABSTRACT
%% ------------------------------------------------------------------
\begin{tcbverbatimwrite}{tmp_\jobname_abstract.tex}
\begin{abstract}
A significant milestone in modern gradient-based optimization was achieved with the development of~\textit{Nesterov's accelerated gradient descent} (\texttt{NAG}) method. This forward-backward technique has been further advanced with the introduction of its proximal generalization, commonly known as the \textit{fast iterative shrinkage-thresholding algorithm} (\texttt{FISTA}), which enjoys widespread application in image science and engineering. Nonetheless, it remains unclear whether both~\texttt{NAG} and~\texttt{FISTA} exhibit linear convergence for strongly convex functions. Remarkably, these algorithms demonstrate convergence without requiring any prior knowledge of strongly convex modulus, and this intriguing characteristic has been acknowledged as an open problem in the comprehensive review~\cite[Appendix B]{chambolle2016introduction}. In this paper, we address this question by utilizing the high-resolution ordinary differential equation (ODE) framework. Expanding upon the established phase-space representation, we emphasize the distinctive approach employed in crafting the Lyapunov function, which involves a dynamically adapting coefficient of kinetic energy that evolves throughout the iterations.  Furthermore, we highlight that the linear convergence of both~\texttt{NAG} and~\texttt{FISTA} is independent of the parameter $r$. Additionally, we demonstrate that the square of the proximal subgradient norm likewise advances towards linear convergence.

%For modern gradient-based optimization, a developmental landmark is Nesterov's accelerated gradient descent method, which is proposed in~\cite{nesterov1983method}, so shortened as~\texttt{Nesterov-1983}. Afterward, one of the important progresses is its proximal generalization, named the fast iterative shrinkage-thresholding algorithm (\texttt{FISTA}), which is widely used in image science and engineering. However, it is unknown whether both~\texttt{Nesterov-1983} and~\texttt{FISTA} converge linearly on the strongly convex function, which has been listed as the open problem in the comprehensive review~\cite[Appendix B]{chambolle2016introduction}. In this paper, we answer this question by the use of the high-resolution differential equation framework.  Along with the phase-space representation previously adopted, the key difference here in constructing the Lyapunov function is that the coefficient of the kinetic energy varies with the iteration. Furthermore, we point out that the linear convergence of both the two algorithms above has no dependence on the parameter $r$ on the strongly convex function. Meanwhile, it is also obtained that the proximal subgradient norm converges linearly.
\end{abstract}

\begin{keywords}
\texttt{NAG}, \texttt{FISTA}, $\mu$-strongly convex function, Lypunov function, phase-space representation
\end{keywords}

\end{tcbverbatimwrite}
\begin{abstract}
A significant milestone in modern gradient-based optimization was achieved with the development of~\textit{Nesterov's accelerated gradient descent} (\texttt{NAG}) method. This forward-backward technique has been further advanced with the introduction of its proximal generalization, commonly known as the \textit{fast iterative shrinkage-thresholding algorithm} (\texttt{FISTA}), which enjoys widespread application in image science and engineering. Nonetheless, it remains unclear whether both~\texttt{NAG} and~\texttt{FISTA} exhibit linear convergence for strongly convex functions. Remarkably, these algorithms demonstrate convergence without requiring any prior knowledge of strongly convex modulus, and this intriguing characteristic has been acknowledged as an open problem in the comprehensive review~\cite[Appendix B]{chambolle2016introduction}. In this paper, we address this question by utilizing the high-resolution ordinary differential equation (ODE) framework. Expanding upon the established phase-space representation, we emphasize the distinctive approach employed in crafting the Lyapunov function, which involves a dynamically adapting coefficient of kinetic energy that evolves throughout the iterations.  Furthermore, we highlight that the linear convergence of both~\texttt{NAG} and~\texttt{FISTA} is independent of the parameter $r$. Additionally, we demonstrate that the square of the proximal subgradient norm likewise advances towards linear convergence.

%For modern gradient-based optimization, a developmental landmark is Nesterov's accelerated gradient descent method, which is proposed in~\cite{nesterov1983method}, so shortened as~\texttt{Nesterov-1983}. Afterward, one of the important progresses is its proximal generalization, named the fast iterative shrinkage-thresholding algorithm (\texttt{FISTA}), which is widely used in image science and engineering. However, it is unknown whether both~\texttt{Nesterov-1983} and~\texttt{FISTA} converge linearly on the strongly convex function, which has been listed as the open problem in the comprehensive review~\cite[Appendix B]{chambolle2016introduction}. In this paper, we answer this question by the use of the high-resolution differential equation framework.  Along with the phase-space representation previously adopted, the key difference here in constructing the Lyapunov function is that the coefficient of the kinetic energy varies with the iteration. Furthermore, we point out that the linear convergence of both the two algorithms above has no dependence on the parameter $r$ on the strongly convex function. Meanwhile, it is also obtained that the proximal subgradient norm converges linearly.
\end{abstract}

\begin{keywords}
\texttt{NAG}, \texttt{FISTA}, $\mu$-strongly convex function, Lypunov function, phase-space representation
\end{keywords}


%% ------------------------------------------------------------------
%% END HEADER
%% ------------------------------------------------------------------

\section{Introduction}
\label{sec: intro}

Since the advent of the new century, statistical machine learning has experienced a rapid expansion, establishing itself as a cornerstone for the development of optimization algorithms. One of the fundamental challenges in this arena is the taks of minimizing an objective function without the imposition of constraints, a scenario typically referred to as an unconstrained optimization problem. This problem is mathematically expressed as 
\[
\min_{x \in \mathbb{R}^d}f(x). 
\]
Gradient-based techniques have risen to prominence due to their computational efficiency and minimal storage requirements, , thereby becoming the method of choice for cutting-edge progress in the field.

When we look back at the historical progression of gradient-based methods, a notable landmark that stands out is~\textit{Nesterov's accelerated gradient descent} (\texttt{NAG}) method. Starting with any initial point $y_0 = x_0 \in \mathbb{R}^d$,~\texttt{NAG} proceeds through the following iterative scheme:
\[
\left\{ \begin{aligned}
        & x_{k} = y_{k-1} - s\nabla f(y_{k-1}), \\
        & y_{k} = x_{k} + \frac{k-1}{k+r}(x_k - x_{k-1}),
        \end{aligned} \right.
\]
 where $s>0$ denotes the fixed step size. Originally proposed in~\cite{nesterov1983method},~\texttt{NAG} is distinguished by its superior convergence performance, particularly regarding convex functions, where it delivers an accelerated convergence rate of $O(1/k^2)$ that represents a substantial improvement over the $O(1/k)$ rate observed in classical gradient descent methods. The underlying mechanism behind the acceleration phenomenon was elucidated with the advent of the high-resolution ordinary differential equation (ODE) framework, as proposed in~\cite{shi2022understanding}. This innovative framework employs phase-space representation in conjunction with Lyapunov analysis to decode the enhanced convergence rates applicable to both the function value and the square of the gradient norm. Subsequent refinements, such as those described in~\cite{chen2022gradient}, simplify this theoretical framework into a concise step while offering an alternative verification through the concept of implicit-velocity phase-space representation. Furthermore, the discovery of a refined proximal inequality, arising from a key observation, paved the way for generalizing the high-resolution ODE framework to encompass composite optimization problems. These advancements have ultimately led to the development of the~\textit{fast iterative shrinkage-thresholding algorithm}~(\texttt{FISTA}), expounded in~\cite{li2022proximal}. Additionally, with a small modification, the accelerated convergence rate for both~\texttt{NAG} and~\texttt{FISTA} has been applied to the iterates, as further explored in~\cite{chambolle2015convergence}. 
 
%  In, a simplified derivation reduced the computation to a single step, offering an alternative proof from the perspective of implicit-velocity phase-space representation. 
% the question of whether its proximal generalization,~\texttt{FISTA}, can achieve a linear convergence rate on strongly convex functions has been listed as an open question

% However, it is worth noting that this variant requires pre-estimation of certain parameters, which can make it impractical to implement in real-world scenarios. 

% The derivation has been greatly simplified to only one-step computation in~\cite{chen2022gradient}, which also provides the other proof from the implicit-velocity phase-space representation. Furthermore, according to the discovery of the improved fundamental proximal inequality based on a key observation, a bridge is built to generalize the high-resolution differential equation framework across over the composite optimization in~\cite{li2022proximal}, which is the so-called~\texttt{FISTA}, the fast iterative shrinkage-thresholding algorithm. 

\begin{figure}[htb!]
\centering
\begin{subfigure}[t]{0.495\linewidth}
\centering
\includegraphics[scale=0.16]{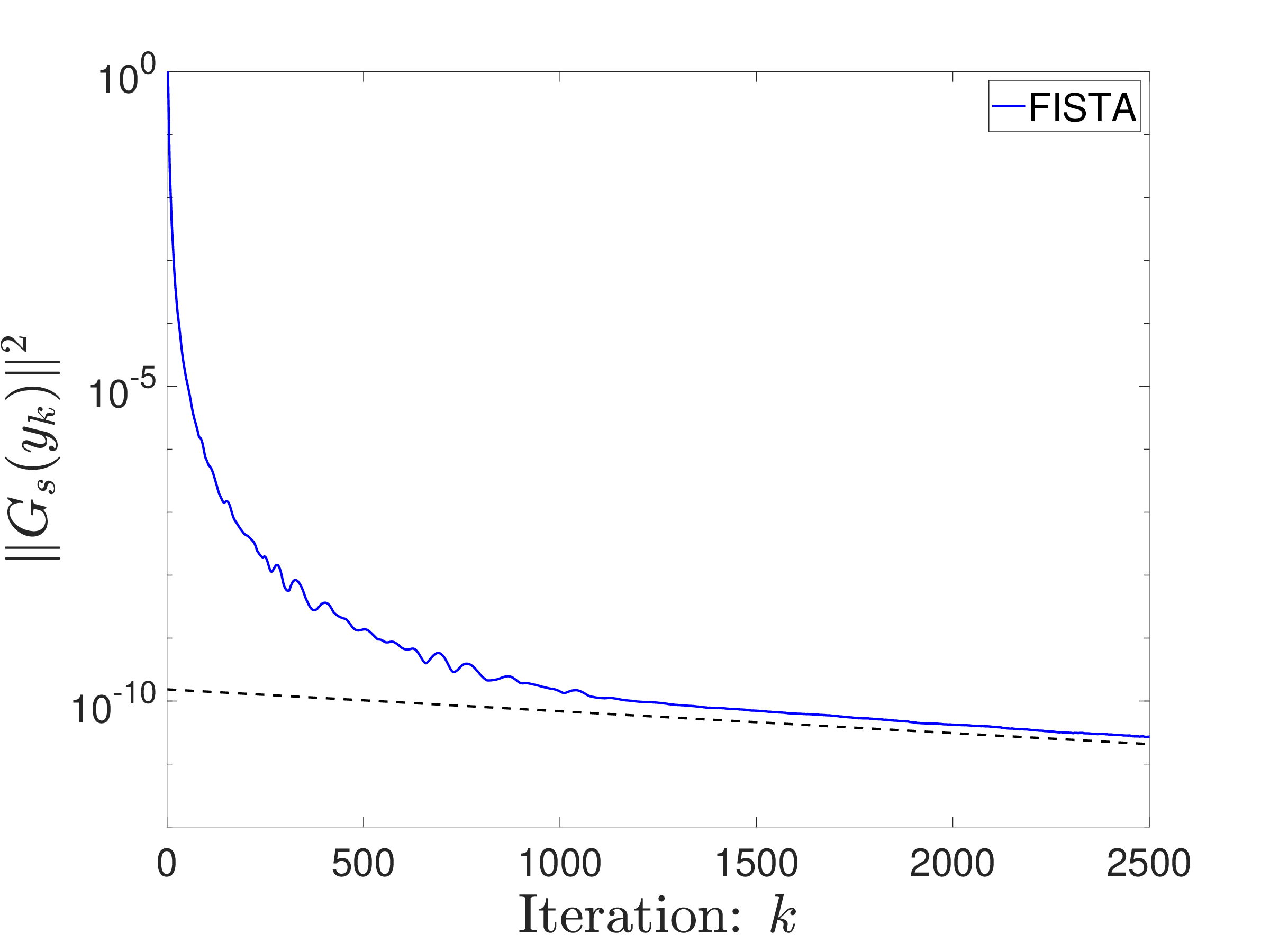}
\caption{Iteration: Start From $k=0$}
\end{subfigure}
\begin{subfigure}[t]{0.495\linewidth}
\centering
\includegraphics[scale=0.16]{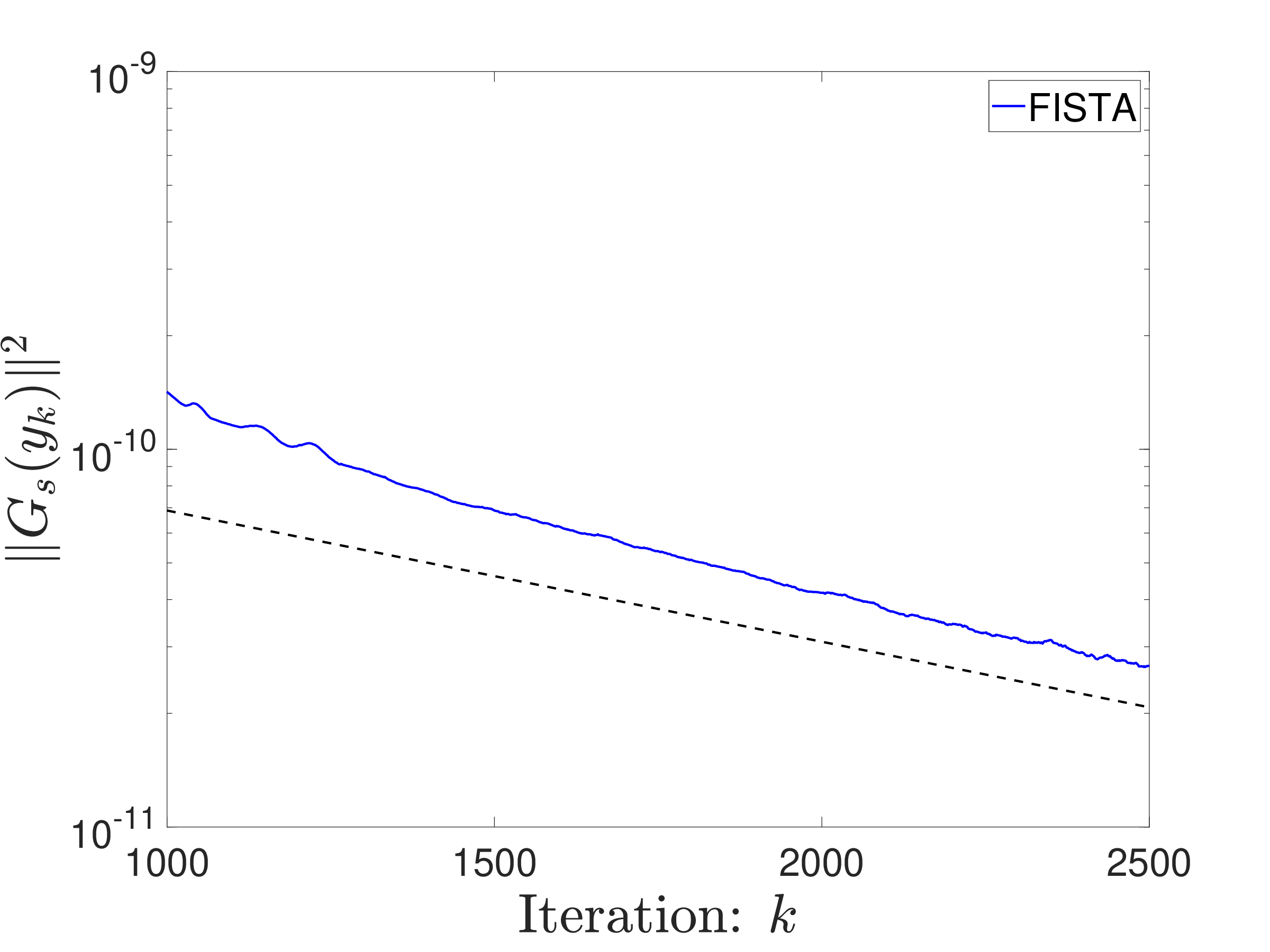}
\caption{Iteration: Start From $k=1000$}
\end{subfigure}
\caption{Iterative progression of the square of the proximal subgradient norm throughout the application of~\texttt{FISTA} for image deblurring, as demonstrated in~\cite{li2022proximal}. } 
\label{fig: fista-strongly}
\end{figure}

In practical scenarios involving convex objective functions, it is crucial to recognize that their associated Hessian matrices are devoid of zero eigenvalues. More commonly, the spectrum of these matrices is characterized by a small ratio between the smallest and largest eigenvalues, a condition that leads to the functions being labeled as ``ill-conditioned'' within the lexicon of the optimization community. To be explicit, for optimal outcomes, the objective function should manifest strong convexity rather than mere convexity in general. A variant of~\texttt{NAG},  proposed in~\cite{nesterov1998introductory},  is designed to fast-track convergence for strongly convex functions. It is noteworthy, however, that this iteration is contingent upon the advanced estimation of certain parameters, a requirement that may impede its practical utility in real-world contexts.  In stark contrast, the original~\texttt{NAG} algorithm operates independently of any foreknowledge regarding the modulus of strong convexity, prompting inquiries into its specific convergence rate when applied to strongly convex functions. These inquiries equally pertain to~\texttt{FISTA}, the proximal generalization of \texttt{NAG}. Whether these algorithms,~\texttt{NAG} and~\texttt{FISTA}, can achieve linear convergence for strongly convex functions is currently recognized as an open question, as outlined in~\cite[Appendix B]{chambolle2016introduction}.  For a more concrete understanding of the convergence rates, we reference the case study presented in~\cite{li2022proximal},  where~\texttt{FISTA} is utilized to deblur an image featuring an elephant.  As a starting point for our analysis, we scrutinize the numerical pattern exhibited in~\Cref{fig: fista-strongly}, which delineates how the square of the proximal subgradient norm varies across successive iterations.

\vspace{0.1cm}

%\paragraph{Main contribution} 

%%%%%%%%%%%%%%%%%%%%%%%%%%%%%%%%%%%%%%%%%%%%%%%%%%%%%%%%%%%%%
\subsection{Overview of contribution}
\label{subsec: contribution} 

In this study, we make use of the high-resolution ODE framework, as previously established in a series of studies~\cite{shi2022understanding, shi2019acceleration, chen2022gradient, chen2022revisiting, chen2023underdamped, li2022linear, li2022proximal}, to address the open question posed in~\cite[Appendix B]{chambolle2016introduction}. Our primary approach revolves around leveraging Lyapunov analysis. The main contributions of our research are outlined below. 

\begin{itemize}
\item We establish the linear convergence of \texttt{NAG} for (smooth) strongly convex functions by formulating an innovative Lyapunov function. What sets our approach apart from previous methods is the incorporation of a dynamically adapting coefficient for the kinetic energy that evolves throughout each iteration, a distinctive feature absent in preceding methods. Importantly, our results indicate that this achieved linear convergence does not depend on the parameter $r$, distinguishing it from the patterns observed in convex scenarios.

% where the key ingredient different from previous ones is that the coefficient of the kinetic energy varies with the iteration. Meanwhile, another difference from the convex case, the linear convergence has no dependence on the parameter $r$ on the strongly convex function.

\item  Furthermore, we refine a key inequality associated with strong convexity to encompass the proximal setting. This enhancement effectively reconciles the fields of smooth and composite optimization, mending a theoretical delineation. With the aid of the implicit-velocity phase-space representation, the inventive Lyapunov function that we devise guarantees the linear convergence of function values within~\texttt{FISTA} and clearly delineates the linear convergence of the square of the proximal subgradient norm. 

\end{itemize}
% has 
%
%
% $\mu$-strongly
%Open question of ISTA and FISTA in~\cite[Appendix B]{chambolle2016introduction}
%
%
%\begin{tcolorbox}
%\begin{itemize}
%\item  Does the~\texttt{Nesterov-1983} converge linearly for the $\mu$-strongly convex function?
%\item Does still work for~\texttt{ISTA} and~\texttt{FISTA}~\cite[Appendix B]{chambolle2016introduction}?
%\end{itemize}
%\end{tcolorbox}
%
%
%
%
%along the iterative trajectory of~\texttt{Nesterov-1983}, the error's convergence rate is improved to $O(1/k^2)$ for the general convex function.\footnote{The error: the difference between the function value  and the minimal value, $f(x_{k}) - f(x^\star)$} According to the existence of the parameter $\mu$ in the momentum coefficient, the variant of~\texttt{Nesterov-1983} accelerating the convergence rate on the $\mu$-strongly convex function\footnote{, $\mu$-strongly convex} does not work in practice. 
%%%%%%%%%%%%%%%%%%%%%%%%%%%%%%%%%%%%%%%%%%%%%%%%%%%%%%%%%%%%%%%%%%%%%%%%%
\subsection{Intuitive analysis on a quadratic function}
\label{subsec: intuitive-analysis}

%Step size $s=1$, Objective function
%\[
%f(x_1,x_2) = 2\times 10^{-2}x_1^2 + 5 \times 10^{-4}x_2^2
%\]

\begin{figure}[htb!]
\centering
\begin{subfigure}[t]{0.495\linewidth}
\centering
\includegraphics[scale=0.16]{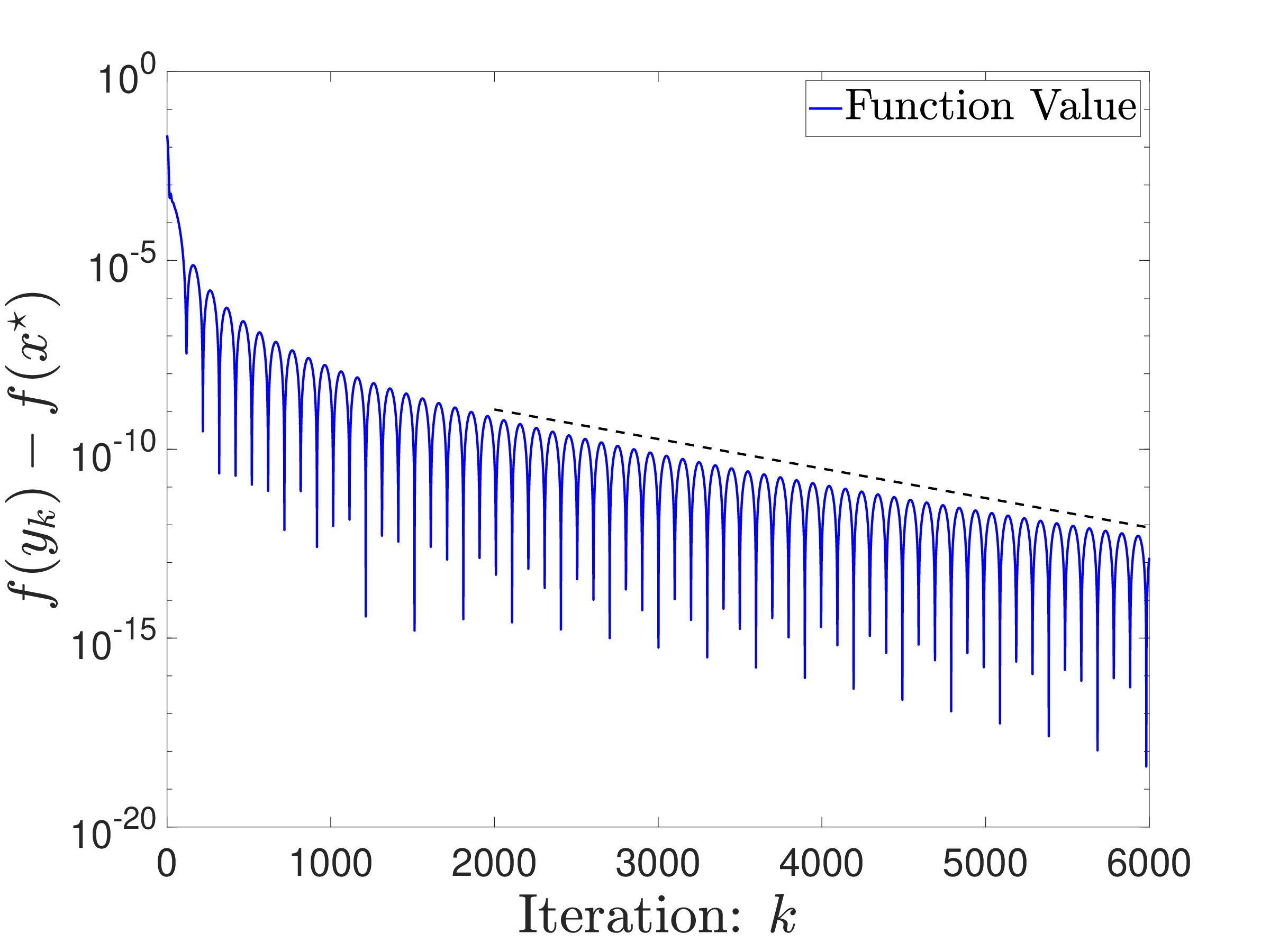}
\caption{Function Value}
\end{subfigure}
\begin{subfigure}[t]{0.495\linewidth}
\centering
\includegraphics[scale=0.16]{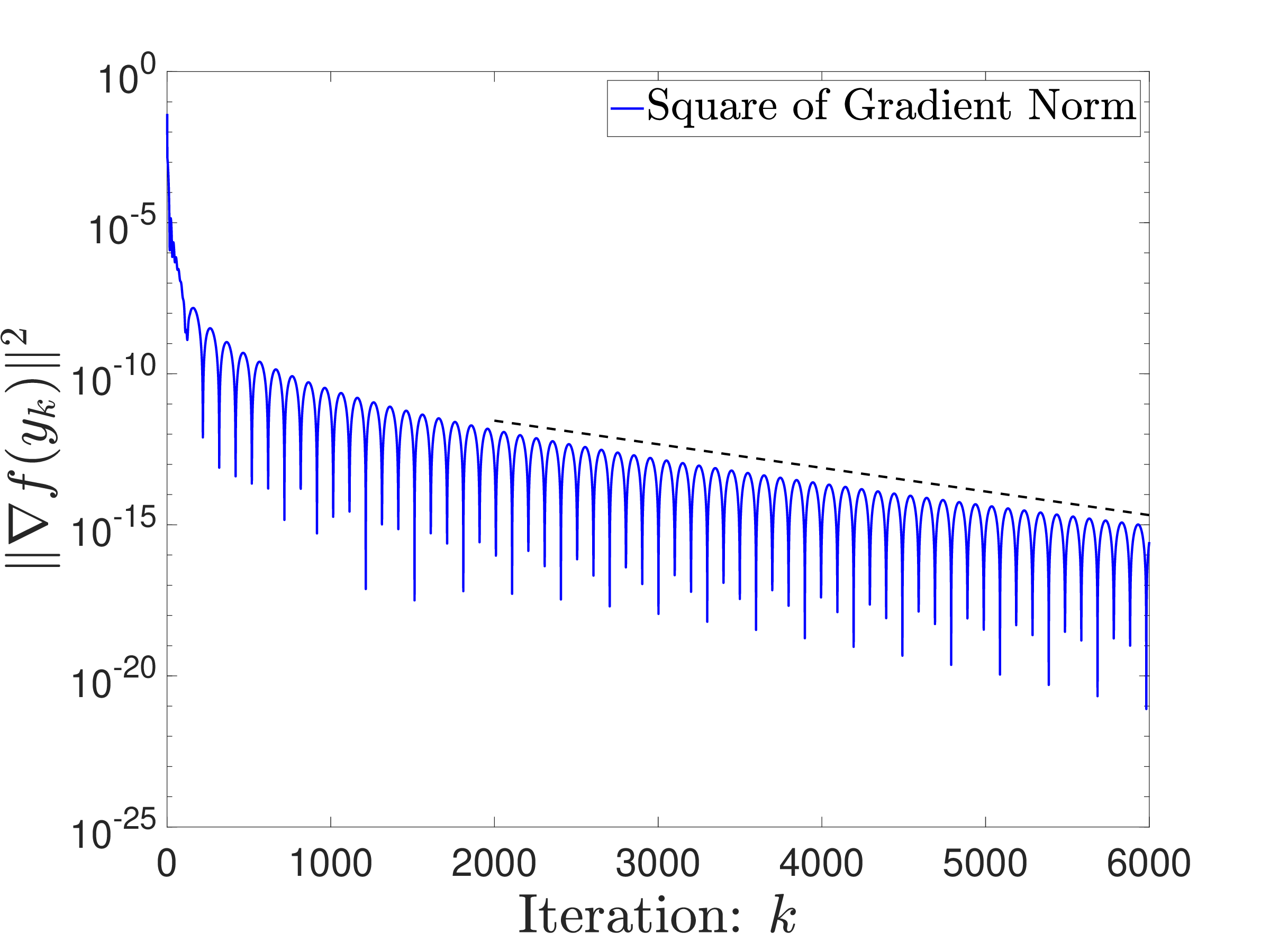}
\caption{Square of Gradient Norm}
\end{subfigure}
\caption{Iterative progression of both the function value and the square of the gradient norm throughout the application of~\texttt{NAG} with a step size of $s=1$, on the quadratic objective function $f(x_1,x_2) = 2\times 10^{-2}x_1^2 + 5 \times 10^{-4}x_2^2$. } 
\label{fig: nag-quadratic}
\end{figure}

To enhance our intuitive understanding, we initiate our discussion with a visual representation as provided in~\Cref{fig: nag-quadratic}. This demonstrative portrayal reveals how both the function value and the square of the gradient norm converge across successive iterations of~\texttt{NAG} when employed on an ill-conditioned quadratic objective function. This figure effectively highlights the central attribute of linear convergence as it pertains to both the function value and the square of the gradient norm.
%\begin{figure}[htbp!]
%\centering
%%\begin{subfigure}[t]{0.48\linewidth}
%%\centering
%\includegraphics[scale=0.28]{fig/nag1.eps}
%%\caption{Error's Convergence}
%%\end{subfigure}
%%\begin{subfigure}[t]{0.48\linewidth}
%%\centering
%%\includegraphics[scale=0.18]{fig/nag2.eps}
%%\caption{Gradient Norm Minimization}
%%\end{subfigure}
%\caption{The error's iterative behavior of~\texttt{Nesterov-1983} with the step size $s=1$ on the quadratic objective function $f(x_1,x_2) = 2\times 10^{-2}x_1^2 + 5 \times 10^{-4}x_2^2$. } 
%\label{fig: nag-quadratic}
%\end{figure}
 
%Since there is no interaction along the eigenvectors' directions of the Hessian for the quadratic function,

For a quadratic function, its Hessian is invariant with respect to the position $x$, leading to consistent eigendirections that are mutually orthogonal. For the sake of simplicity, let us consider the objective function to be one-dimensional, taking the general form $f(x) = \frac12 \mu x^2$, where $\mu > 0$.  By incorporating this general form into~\texttt{NAG}, the iteration scheme is postulated as
\[
\left\{ \begin{aligned}
        & x_{k} = y_{k-1} - \mu s y_{k-1}, \\
        & y_{k} = x_{k} + \frac{k-1}{k+r}(x_k - x_{k-1}).
        \end{aligned} \right.
\]
By performing some elementary operations of linear transformations, we can recast the iteration in vector form as 
\[
\begin{pmatrix}
x_k \\ y_k
\end{pmatrix}  =  \begin{pmatrix}
0 & 1 - \mu s \\ - \frac{k-1}{k+r} & \frac{2k+r-1}{k+r} \cdot (1 - \mu s)
\end{pmatrix}
\begin{pmatrix}
x_{k-1} \\ y_{k-1}
\end{pmatrix},
\]
where the eigenvalues of the iterative matrix adhere to the characteristic quadratic polynomial equation: 
\[
\lambda^2 -  \frac{2k+r-1}{k+r} (1 - \mu s) \cdot \lambda + \frac{k-1}{k+r}(1 - \mu s) = 0. 
\]
%\end{equation}
As the iterations progress proliferate infinity ($k \rightarrow \infty$), the discriminant of this quadratic equation asymptotocally converges to
%\[
%\lambda^2 -  \frac{2k+r-1}{k+r} (1 - s\theta) \cdot \lambda + \frac{k-1}{k+r}(1 - s\theta) = 0
%\]
\[
 \left(\frac{2k+r-1}{k+r}\right)^2 (1 - \mu s)^2 - 4 \left(\frac{k-1}{k+r}\right)(1 - \mu s) \rightarrow -4\mu s (1 - \mu s). 
\]
Given $\mu s \in (0,1)$, the asymptotic expressions for the roots of this quadratic equation become:
\[
\lambda_{1,2} \rightarrow 1 - \mu s \pm \sqrt{\mu s(1 - \mu s)}i \qquad \text{as} \quad k \rightarrow \infty.
\]
This demonstrates that~\texttt{NAG} consistently achieves linear convergence at an asymptotic rate of $\sqrt{1 - \mu s}$,  analogous to the convergence rate exhibited in classic gradient descent, as explicated in~\cite[Theorem A.4]{shi2019acceleration}.

\begin{figure}[htb!]
\centering
\begin{subfigure}[t]{0.495\linewidth}
\centering
\includegraphics[scale=0.16]{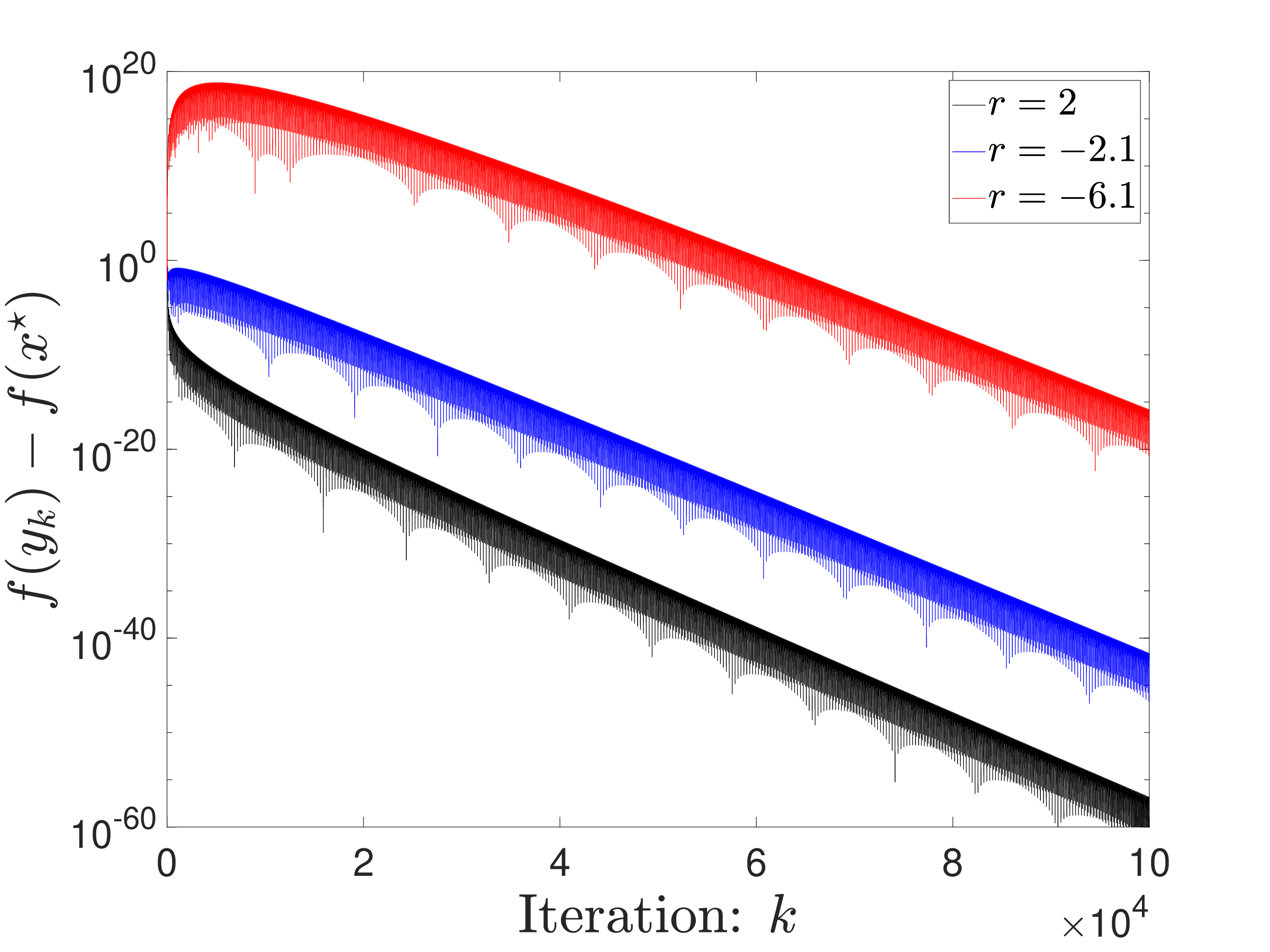}
\caption{Function Value}
\end{subfigure}
\begin{subfigure}[t]{0.495\linewidth}
\centering
\includegraphics[scale=0.16]{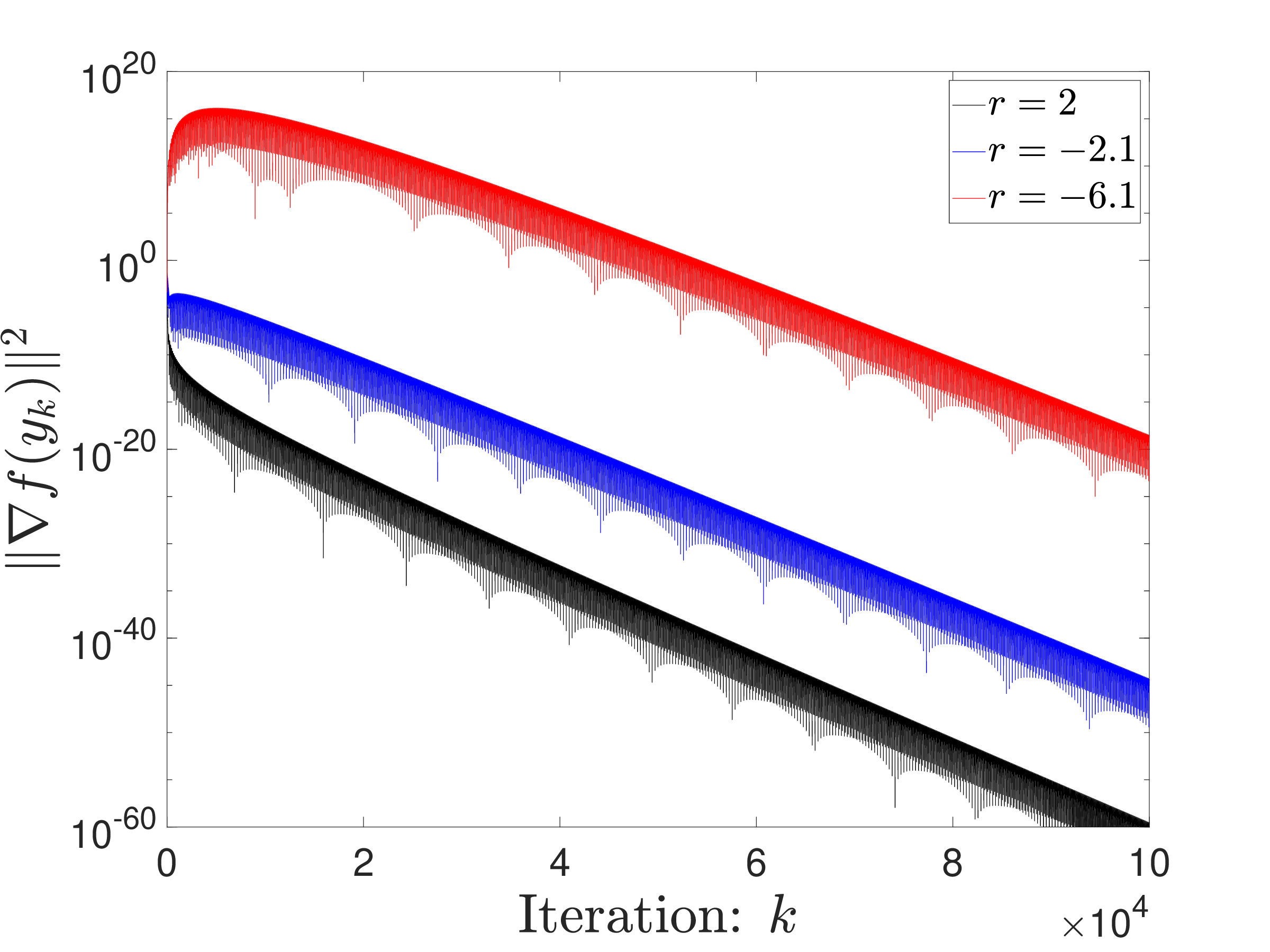}
\caption{Square of Gradient Norm}
\end{subfigure}
\caption{Illustration of linear convergence of~\texttt{NAG}, demonstrated by tracking both the function value and the square of the gradient norm, as it varies with the parameter $r$ in the same scenario as depicted in~\Cref{fig: nag-quadratic}.   } 
\label{fig: r-negtive}
\end{figure}

Furthermore, it becomes apparent that the linear convergence rate of~\texttt{NAG}  is rarely influenced by the choice of the parameter $r$.  This observation is supported by the results of a numerical experiment in which the parameter $r$ is assigned across three distinct values, with the outcomes depicted in~\Cref{fig: r-negtive}. When the parameter $r$ is positive, it is observed that both the function value and the square of the gradient norm always exhibit convergence. In contrast, when the parameter $r$ is negative, there is an initial uptick in both the function value and the square of the gradient norm before they shift toward eventual convergence. As demonstrated in~\Cref{fig: r-negtive}, after exceeding $20,000$ iterations, the linear convergence rate retains comparable characteristics across diverse settings of the parameter $r$, regardless of whether assessing the function value or the square of the gradient norm. While distinct patterns of iteration may manifest in the early stages, a consistent rate of linear convergence is eventually achieved as the iterations progress. To preclude complications associated with a zero denominator in the momentum coefficient, which may occur when $k + r$ equals zero, non-integer values are deliberately selected for the negative parameter $r$ in the numerical experiment illustrated in~\Cref{fig: r-negtive}. Traditional numerical approaches usually begin with the initial iteration set to zero; nonetheless, to bypass the zero denominator concern, the experiment may alternatively begin at the $(-r+1)$-th iteration, corresponding to $k=-r+1$. Hence, by adopting this approach, the parameter $r$ can be retained as an integer value for negative instances as well, thereby preventing any operational issues.

% is set not integer in order to $k + r \neq 0$, if we set $r$ integer, start from $r+1$ 

%\[
%\begin{pmatrix}
%1 & 0 \\ - \frac{2k+r-1}{k+r} & 1
%\end{pmatrix}
%\begin{pmatrix}
%x_k \\ y_k
%\end{pmatrix} = \begin{pmatrix}
%0 & 1 - s\theta \\ - \frac{k-1}{k+r} & 0
%\end{pmatrix}
%\begin{pmatrix}
%x_{k-1} \\ y_{k-1}
%\end{pmatrix} 
%\]
%
%\[
%\begin{pmatrix}
%x_k \\ y_k
%\end{pmatrix} = 
%\begin{pmatrix}
%1 & 0 \\ \frac{2k+r-1}{k+r} & 1
%\end{pmatrix}
%\begin{pmatrix}
%0 & 1 - s\theta \\ - \frac{k-1}{k+r} & 0
%\end{pmatrix}
%\begin{pmatrix}
%x_{k-1} \\ y_{k-1}
%\end{pmatrix} =  \begin{pmatrix}
%0 & 1 - s\theta \\ - \frac{k-1}{k+r} & \frac{2k+r-1}{k+r} \cdot (1 - s\theta)
%\end{pmatrix}
%\begin{pmatrix}
%x_{k-1} \\ y_{k-1}
%\end{pmatrix}
%\]

%%%%%%%%%%%%%%%%%%%%%%%%%%%%%%%%%%%%%%%%%%%%%%%%%%%%%%%%%%%%%
\subsection{Related works and organization}
\label{subsec: related-work} 

The history of first-order optimization algorithms, specifically focusing on acceleration, can be traced back to the inception of the Ravine method. This method, a two-step gradient strategy, was designed to outperform the classical single-step gradient descent~\cite{gelfand1962some}. The renowned~\texttt{NAG} method and its variant were introduced in~\cite{nesterov1983method} and~\cite{nesterov1998introductory}, catalyzing a revolution in the phenomena of global acceleration within convex optimization. The discovery of gradient correction, unveiling the underlying mechanism of acceleration, was recognized through the comparison of the~\texttt{NAG} method with Polyak's heavy-ball method, within the framework of high-resolution ODEs proposed in~\cite{shi2022understanding}. It was further identified in~\cite{chen2022gradient, chen2022revisiting} that this high-resolution ODE framework is specially tailored for~\texttt{NAG} and its variant, elucidating how gradient correction is effectively realized through implicit velocity. Advances in the fundamental proximal inequality paved the way for the expansion of this high-resolution ODE framework to address composite optimization, applicable to both convex and strongly convex functions, as reported in~\cite{li2022proximal, li2022linear}. In addition, the completion of this tailor-made framework to accommodate the underdamped case where $r<2$ was achieved in~\cite{chen2023underdamped}.

% Until the high-resolution differential equation framework proposed in~\cite{shi2022understanding}, the veil behind the acceleration phenomenon is lifted by the discovery of gradient correction by the comparison with Polyak's heavy-ball method in~\cite{polyak1964some}. Recently, the high-resolution differential equation framework is discovered further to be tailor-made for~\texttt{NAG} and its variant in~\cite{chen2022gradient, chen2022revisiting}, which also finds that the other form to manifest the effect of gradient correction is the implicit velocity. Based on improving the fundamental proximal inequality, the high-resolution differential equation framework is generalized to the composite optimization for both the general and strongly convex functions in~\cite{li2022proximal, li2022linear}. In addition, this tailor-made framework is completed in~\cite{chen2023underdamped} for the underdamped case $r<2$. 
 
Over the past decade, it has attracted considerable attention towards research utilizing ODEs with Hessian-driven damping to explore forward-backward algorithms, with early works published in~\cite{attouch2012second, attouch2014dynamical}. The derivation of a low-resolution ODE as a model for~\texttt{NAG} came from employing techniques borrowed from numerical analysis in~\cite{su2016differential}. This laid the groundwork for a multitude of research pathways, including studies on the faster convergence rate of function values~\cite{attouch2016rate}, the fusion of the low-resolution ODE and Hessian-driven damping to develop continuous dynamics~\cite{attouch2016fast}, the use of variational method~\cite{wibisono2016variational}, applications of Lyapunov analysis~\cite{wilson2021lyapunov}, and investigations into the implicit-velocity perspective~\cite{muehlebach2019dynamical}. The unveiling of the underlying mechanisms was achieved with the advent of high-resolution numerical approximation techniques, as outlined in~\cite{shi2022understanding}. Utilizing phase-space representation and Lyapunov analysis, this work bridged the gap between continuous dynamics and discrete algorithms, casting light on the interplay between the two areas. Additionally, several studies have inspected the acceleration phenomenon with an eye on numerical stability, as detailed in~\cite{luo2022differential, zhang2021revisiting}. Tthe momentum-based scheme has also been extrapolated to the domain of residual neural networks~\cite{xia2021heavy}, further expanding its applicability.

The remainder of this paper is organized as follows. In~\Cref{sec: prelim}, we set forth some basic notations and definitions to serve as preliminaries. In~\Cref{sec: high-resolution}, we offer a proof predicated on the gradient-correction high-resolution ODE that addresses the continuous dynamics.  In~\Cref{sec: linear-convergence}, we are dedicated to proving the linear convergence of~\texttt{NAG} for strongly convex objective functions, employing the iteration-varying Lyapunov function. Finally, in~\Cref{sec: conclusion-further-works}, we discuss the distinctions between the proofs for the continuous ODE and the discrete algorithm and propose potential avenues for future research.

\section{Preliminaries}
\label{sec: prelim} 

The notations employed within this paper mostly aligns with those found in~\cite{nesterov1998introductory, shi2022understanding, li2022linear}, with slight modifications tailored to our context. Let $\mathcal{F}^0(\mathbb{R}^d)$ denote the class of continuous convex functions on $\mathbb{R}^d$; that is, $g \in \mathcal{F}^{0}$ if it fulfills the inequality
\[
g\left( \alpha x + (1 - \alpha)y \right) \leq \alpha g(x) + (1 - \alpha)g(y)
\] 
for any $x,y \in \mathbb{R}^d$ and $\alpha \in [0, 1]$. The subclass $\mathcal{F}^1_L(\mathbb{R}^d) \subseteq \mathcal{F}^0(\mathbb{R}^d)$ consists of functions whose gradients are well-defined everywhere and adheres to the global Lipschitz condition. Thus, $f \in \mathcal{F}^{1}_{L}$ if $f \in \mathcal{F}^{0}$ and it satisfies
\[
\| \nabla f(x) - \nabla f(y) \| \leq L \| x - y \|
\] 
for any $x, y \in \mathbb{R}^d$.\footnote{Throughout this paper, the notation $\|\cdot\|$ specifically refers to the $\ell_2$-norm or Euclidean norm, denoted as $\|\cdot\|_2$. It is worth noting that the subscript $2$ is often omitted for convenience unless otherwise noted.} We also denote $\mathcal{S}_{\mu,L}^{1}(\mathbb{R}^d)$ as the subclass of $\mathcal{F}^1_L(\mathbb{R}^d)$ with each member being $\mu$-strongly convex for some $0 < \mu \leq L$. In other words, $f \in \mathcal{S}_{\mu,L}^{1}(\mathbb{R}^d)$ if $f \in \mathcal{F}^1_L(\mathbb{R}^d)$ and it holds that
\[
f(y) \geq f(x) +  \langle \nabla f(y), y - x \rangle + \frac{\mu}{2} \|y - x\|^2
\]
for any $x, y \in \mathbb{R}^d$. Furthermore, $\mathcal{S}_{\mu,L}^{2}(\mathbb{R}^d)$ refers to a subclass of $\mathcal{S}_{\mu,L}^{1}(\mathbb{R}^d)$, which encompasses functions possessing a continuous Hessian. Finally,  the term $x^\star$ is used to denote its unique minimizer of those functions. 

Additionally, we consider a composite function $\Phi= f + g$, where $f \in \mathcal{S}_{\mu,L}^1$ and $g \in \mathcal{F}^0$.  This is analogous to~\cite{beck2009fast, su2016differential}, where the $s$-proximal operator and the $s$-proximal subgradient operator are defined as follows.  
\begin{definition}
\label{defn: proximal-subgradient}
Let the step size satisfy $s \in (0, 1/L)$. For any $f\in \mathcal{S}_{\mu,L}^1$ and $g \in \mathcal{F}^0$, the $s$-proximal operator is defined as
\begin{equation}
\label{eqn: proximal-operator}
        P_s(x) := \mathop{\arg\min}_{y\in\mathbb{R}^d}\left\{ \frac{1}{2s}\left\| y - \left(x - s\nabla f(x)\right) \right\|^2 + g(y) \right\},
\end{equation}
for any $x \in \mathbb{R}^d$. Furthermore, the $s$-proximal subgradient operator is defined as
\begin{equation}
\label{eqn: subgradient-operator}
G_s(x): = \frac{x - P_s(x)}{s}
\end{equation}
for any $x \in \mathbb{R}^d$.
\end{definition}

When $g$ simplifies to the $\ell_1$-norm, i.e., $g(x) = \lambda \|x\|_1$,\footnote{The $\ell_1$ norm, denoted as $\|\cdot\|_1$, is defined as
\[ \|x\|_{1} = \sum_{i=1}^{d} |x_i|,\]
for any $x \in \mathbb{R}^d$.} we can derive the closed-form silution for the $s$-proximal operator~\eqref{eqn: proximal-operator} at any $x \in \mathbb{R}^d$ for the particular instance as
\[
P_s(x)_i = \big(\left|\left(x - s\nabla f(x)\right)_i\right| - \lambda s\big)_+ \text{sgn}\big(\left(x - s\nabla f(x)\right)_i\big),
\]
where $i=1,\ldots,d$ represents the index.

\section{The gradient-correction high-resolution ODE}
\label{sec: high-resolution}

%We demonstrate the linear convergence of~\texttt{NAG} by use of its implicit-velocity phase-space representation in the last section.

In this section, we delve into the area of continuous convergence rates. Our primary focus lies on the gradient-correction high-resolution ODE, which is articulated in~\cite{shi2022understanding}. Serving as a continuous analog of~\texttt{NAG}, the gradient-correction high-resolution ODE  is expressed as:
\begin{equation}
\label{eqn: grad-correction-ode}
    \ddot{X} + \dfrac{3}{t}\dot{X} + \sqrt{s}\nabla^2 f(X) \dot{X} + \brak{1 + \dfrac{3\sqrt{s}}{2t}}\nabla f(X) = 0,
\end{equation}
with any initial $(X(0), \dot{X}(0)) = (x_{0},v_{0}) \in \mathbb{R}^d \times \mathbb{R}^d$. Given that $f \in \mathcal{S}_{\mu,L}^{2}$, it is certain that  the eigenvalue of the Hessian is always greater than or equal to~$\mu$, which is denoted as $\lambda(\nabla^2 f(x))  \geq \mu$. This observation suggests that the damping term remains consistently substantial, thereby hinting at the potential for linear convergence.   To determine the exact convergence rate, we turn to the principled approach of constructing a Lyapunov function, as demonstrated in~\cite{chen2022revisiting}.
\begin{itemize}
\item[(\textbf{I})] We initiate our analysis with the mixed energy, inspired by the convex case mentioned in~\cite[(4.36)]{shi2022understanding}, expressed as
\begin{equation}
\label{eqn: con-mix}
\E_{\textbf{mix}} = \dfrac{1}{2}\|t\dot{X} + 2(X - x^\star) + t\sqrt{s}\nabla f(X)\|^2.
\end{equation}
By using the gradient-correction high-resolution ODE~\eqref{eqn: grad-correction-ode}, we can calculate its time derivative as
\begin{align}
        \frac{d \E_{\textbf{mix}}}{d t} & = \inner{t\dot{X} + 2(X - x^\star) + t\sqrt{s}\nabla f(X), -\brak{t + \frac{\sqrt{s}}{2}}\nabla f(X)} \nonumber \\
                                            & = \underbrace{-t\brak{t+\frac{\sqrt{s}}{2}} \big\langle\dot{X}, \nabla f(X)\big\rangle}_{\mathbf{I}} - (2t + \sqrt{s})\big\langle\nabla f(X), X - x^\star\big\rangle \nonumber \\
                                            & \mathrel{\phantom{=}}- \sqrt{s}t\brak{t + \frac{\sqrt{s}}{2}}\norm{\nabla f(X)}^2. \label{eqn: time-derive-mix}
\end{align}

\item[(\textbf{II})] Adhering to the principled approach of constructing a Lyapunov function as described in~\cite{chen2022revisiting}, we consider the kinetic energy. Different from~\cite[(2.4)]{shi2022understanding}, here the kinetic energy includes a time-varying coefficient as
\begin{equation}
\label{eqn: con-kin}
\E_{\textbf{kin}} = \dfrac{\tau(t)}{2}\|\dot{X}\|^2.
\end{equation}
Taking into account  the gradient-correction high-resolution ODE~\eqref{eqn: grad-correction-ode}, we can calculate the time derivative of the kinetic energy as
\begin{align}
        \frac{d \E_{\textbf{kin}}}{d t} & = \tau(t) \big\langle\ddot{X}, \dot{X}\big\rangle + \frac{\dot{\tau}(t)}{2}\|\dot{X}\|^2 \nonumber\\
                                                    & = -\left( \frac{3\tau(t)}{t} - \frac{\dot{\tau}(t)}{2} \right)\|\dot{X}\|^2 - \sqrt{s}\tau(t)\dot{X}^{\mathrm{T}} \nabla^2 f(X)\dot{X} \nonumber \\ & \mathrel{\phantom{=}}\underbrace{- \tau(t)\brak{1+\frac{3\sqrt{s}}{2t}}\big\langle\dot{X}, \nabla f(X)\big\rangle}_{\mathbf{II}}.\label{eqn: time-derive-kin}
\end{align}
To simplify the coefficient of $\|\dot{X}\|^2$, a good choice is to let $\tau(t)$ to be a power of $t$, such as $\tau(t) = t^{\alpha}$. In order to combine $\mathbf{I}$ and $\mathbf{II}$, it is appropriate for us to choose $\alpha=2$.

\item[(\textbf{III})] Following the principled approach, to eliminate the term involving $\big\langle\dot{X}, \nabla f(X)\big\rangle$, the coefficient of the potential energy should be set accordingly. Specifically, amalgamating terms $\mathbf{I}$ and $\mathbf{II}$ yields the expression $\mathbf{I} + \mathbf{II} = - 2t(t + \sqrt{s}) \big\langle\dot{X}, \nabla f(X) \big\rangle$.  Consequently, the potential energy is constructed as
\begin{equation}
\label{eqn: con-pot}
\E_{\textbf{pot}} = 2t(t + \sqrt{s})\brak{f(X) - f(x^\star)}
\end{equation}
Proceeding with the exploration of the gradient-correction high-resolution ODE~\eqref{eqn: grad-correction-ode}, we calculate its time derivative as follows:
\begin{equation}
\label{eqn: time-derive-pot}
        \frac{d\E_{\textbf{pot}}}{d t}  =\brak{4t+2\sqrt{s}}\brak{f(X) - f(x^\star)} + \underbrace{2t(t + \sqrt{s}) \langle\dot{X}, \nabla f(X)\rangle}_{\mathbf{III}},
\end{equation}
where the resultant term $\mathbf{III}$ in~\eqref{eqn: time-derive-pot} ensures that the combination of $\mathbf{I}$, $\mathbf{II}$, and $\mathbf{III}$ satisfies the equality $\mathbf{I} + \mathbf{II} + \mathbf{III}  = 0$.
\end{itemize}

By integrating the mixed energy~\eqref{eqn: con-mix}, the kinetic energy~\eqref{eqn: con-kin} and the potential energy~\eqref{eqn: con-pot}, we arrive at the following Lyapunov function as 
\begin{equation}\label{eqn: continuous_Lyapunov}
    \E = 2t(t+ \sqrt{s})\brak{f(X) - f(x^\star)}+ \dfrac{t^2}{2}\|\dot{X}\|^2+ \dfrac{1}{2}\|t\dot{X} + 2(X - x^\star) + t\sqrt{s}\nabla f(X)\|^2.
\end{equation}
Equipped with the Lyapunov function~\eqref{eqn: continuous_Lyapunov}, we can conclude this section with the following theorem pertaining to the convergence rates of both the function value and the square of the gradient norm, along with its proof. %provide proof for the convergence rate of a continuous high-resolution ODE
%
%
% We show that the high-resolution equation for NAG
%
%has a linear convergence rate under strong convexity. 
%The above theorem shows that the high-resolution equation for NAG has a linear convergence rate $O(\frac{1}{t^2}e^{-\frac{\mu\sqrt{s}}{4}}t)$ to the global minimum as $t$ large enough. Following the standard steps, we first construt the Lyapunov function as 
%
%Note that the main problem here is that when we calculate the derivative, the order of the coefficients decrease. This prevent us from getting linear convergence rate. In order to handle this, we need the following lemma. 
\begin{theorem}\label{thm: equation}
    Let $f \in \mathcal{S}_{\mu,L}^{2}$. For any step size $0 < s < 1/L$, there exists some time $T = T(\mu, s) > 0$ such that the solution $X = X(t)$ to the gradient-correction high-resolution ODE~\eqref{eqn: grad-correction-ode} satisfies
    \begin{equation}
    \left\{ \begin{aligned}
                    &  f(X) - f(x^\star)\leq \frac{\mathcal{E}(T)}{2t(t+\sqrt{s})} e^{-\frac{\mu\sqrt{s}}{4}(t-T)}, \\
                    & \| \nabla f(X) \|^2 \leq \frac{L\mathcal{E}(T)}{t(t+\sqrt{s})} e^{-\frac{\mu\sqrt{s}}{4}(t-T)},
              \end{aligned}\right.      
    \end{equation}
    for any $t \geq T = \frac{4}{\mu\sqrt{s}}$. 
%Here 
%    \begin{multline*}
%        C(T, s) = \left[2T(T+\sqrt{s})\brak{f(X(T)) - f(x^\star)} + \dfrac{T^2}{2}\norm{\dot{X}(T)}^2 \right. \\\left.+ \dfrac{1}{2}\norm{T\dot{X}(T) + 2(X_T - x^\star) + T\sqrt{s}\nabla f(X(T))}^2\right].
%    \end{multline*}
\end{theorem}

%\begin{lemma}[Polyak-Lojasiewicz inequality]\label{lma: PLinequality}
%    Suppose $f(x)$ is differentiable and $\mu$-strongly convex, then it satisfies
%    \begin{equation}
%        \norm{\nabla f(x)}^2 \ge 2\mu\brak{f(x) - f(x^\star)}.
%    \end{equation}
%\end{lemma}
%The proof can be found in many textbook such as ... We are now ready to prove~\Cref{thm: equation}.

\begin{proof}[Proof of~\Cref{thm: equation}]
%    We divide the calculation of derivative $\frac{\dd}{\dd t}\E(t)$ into three parts. First, we calculate $\frac{\dd}{\dd t} \E_{\textbf{mix}}(t)$ as 
%    \begin{align*}
%        \frac{\dd}{\dd t}\E_{\textbf{mix}}(t) =& \inner{t\dot{X} + 2(X - x^\star) + t\sqrt{s}\nabla f(X), -\brak{t + \frac{\sqrt{s}}{2}}\nabla f(X)} \\
%        =& -t\brak{t+\frac{\sqrt{s}}{2}}\inner{\dot{X}, \nabla f(X)} - (2t + \sqrt{s})\inner{\nabla f(X), X - x^\star} - \sqrt{s}t\brak{t + \frac{\sqrt{s}}{2}}\norm{\nabla f(X)}^2.
%    \end{align*}
%    Then the derivative of $\frac{\dd}{\dd t} \E_{\textbf{kin}}(t)$ is calculated by 
%    \begin{align*}
%        \frac{\dd}{\dd t} \E_{\textbf{mix}}(t) =& t^2\inner{\ddot{X}(t), \dot{X}(t)} + t\norm{\dot{X}(t)}^2\\
%        =& -2t\norm{\dot{X}(t)}^2 - \sqrt{s}t^2\dot{X}^{\mathrm{T}} \nabla^2 f(X)\dot{X} - t\brak{t+\frac{3}{2}\sqrt{s}}\inner{\dot{X}, \nabla f(X)}.
%    \end{align*}
%    As $f(x)$ is $\mu$-strongly convex, we have 
%    \begin{equation*}
%        \dot{X}^{\mathrm{T}} \nabla^2 f(X)\dot{X} \ge \mu \norm{\dot{X}}^2.
%    \end{equation*}
%    Hence 
%    \begin{equation*}
%        \frac{\dd}{\dd t} \E_{\textbf{mix}}(t) \le -\brak{\mu\sqrt{s}t^2 + 2t}\norm{\dot{X}}^2 - t\brak{t + \frac{3}{2}\sqrt{s}}\inner{\dot{X},\nabla f(X)}.
%    \end{equation*}
%    Finally we compute the derivative $\frac{\dd}{\dd t} \E_{\textbf{pot}}(t)$.
%    \begin{align*}
%        \frac{\dd}{\dd t} \E_{\textbf{pot}}(t) = \brak{4t + 2\sqrt{s}}\brak{f(X) - f(x^\star)} + 2t\brak{t + \sqrt{s}}\inner{\dot{X}, \nabla f(X)}.
%    \end{align*}
To prove the theorem, we begin with the Lyapunov function~\eqref{eqn: continuous_Lyapunov}. Its time derivative can be calculated by summing up~\eqref{eqn: time-derive-mix},~\eqref{eqn: time-derive-kin} and~\eqref{eqn: time-derive-pot} as
    \begin{align}
        \frac{d \E}{d t} \le & \brak{4t + 2\sqrt{s}}\brak{f(X) - f(x^\star)} - (2t + \sqrt{s})\inner{\nabla f(X), X - x^\star} \nonumber \\ 
                                    & - 2t\|\dot{X}\|^2 - \sqrt{s}t^2\dot{X}^{\mathrm{T}} \nabla^2 f(X)\dot{X} - \sqrt{s}t\brak{t + \frac{\sqrt{s}}{2}}\norm{\nabla f(X)}^2. \label{eqn: con-total-1}
    \end{align}
    For any $f \in \mathcal{S}_{\mu,L}^{2}$,  it satisfies the $\mu$-strongly convex inequality as
    \begin{equation}
    \label{eqn: strong-cvx-inequ}
    \inner{\nabla f(X), X - x^\star} \geq f(X) - f(x^\star) + \dfrac{\mu}{2}\norm{X - x^\star}^2.
    \end{equation}
   By substituting the $\mu$-strongly convex inequality~\eqref{eqn: strong-cvx-inequ} into the earlier time derivative~\eqref{eqn: con-total-1}, we obtain
     \begin{align}
      \frac{d \E}{d t} \leq & \brak{2t + \sqrt{s}}\brak{f(X) - f(x^\star)}  - \frac{\mu\brak{2t + \sqrt{s}}}{2}\norm{X - x^\star}^2 \nonumber \\ 
                                    & - 2t\|\dot{X}\|^2  - \sqrt{s}t^2\dot{X}^{\mathrm{T}} \nabla^2 f(X)\dot{X} - \sqrt{s}t\brak{t + \frac{\sqrt{s}}{2}}\norm{\nabla f(X)}^2. \label{eqn: con-total-2}
    \end{align}    
   To establish proportionality between the corresponding terms, we can estimate the Lyapunov function using Cauchy-Schwarz inequality as
    \begin{equation}
    \label{eqn: con-lyapunov-estimate}
        \E \leq 2t\brak{t + \sqrt{s}}\brak{f(X) - f(x^\star)} + 2t^2\|\dot{X}\|^2 + 6\norm{X - x^\star}^2 + \frac{3}{2}st^2\norm{\nabla f(X)}^2.
    \end{equation}
For any $f \in \mathcal{S}_{\mu,L}^2$, we have two additional $\mu$-strongly convex inequalities as:
    \begin{subequations}
    \begin{align}
    & \dot{X}^{\mathrm{T}} \nabla^2 f(X)\dot{X} \geq \mu \|\dot{X}\|^2, \label{eqn: strong-cvx-inequ1} \\
    & \norm{\nabla f(X)}^2 \geq 2\mu\brak{f(X) - f(x^\star)}. \label{eqn: strong-cvx-inequ2}
    \end{align}     
    \end{subequations}
    By substituting~\eqref{eqn: strong-cvx-inequ1} and~\eqref{eqn: strong-cvx-inequ2} into~\eqref{eqn: con-total-2}, we can obtain the time derivative as
    \begin{align}
        \frac{d\E}{d t} &\leq - \mu \sqrt{s}\left(t-\frac{1}{\mu\sqrt{s}}\right)\brak{t + \sqrt{s} }\brak{f(X) - f(x^\star)} - \mu\sqrt{s}t^2 \|\dot{X}\|^2 \nonumber \\ &\mathrel{\phantom{\leq}} - \mu t \norm{X - x^\star}^2 - \frac{\sqrt{s}t}{2}\left(t-\frac{1}{\mu\sqrt{s}}\right)\norm{\nabla f(X)}^2 \nonumber \\
                               & \leq - \frac{3\mu t\brak{\sqrt{s}t + s }}{4}\brak{f(X) - f(x^\star)} - \mu\sqrt{s}t^2 \|\dot{X}\|^2 \nonumber \\
                               & \mathrel{\phantom{\leq}} - \frac{4}{\sqrt{s}} \norm{X - x^\star}^2 - \frac{3\sqrt{s}t^2}{8}\norm{\nabla f(X)}^2,  \label{eqn: con-total-3}
    \end{align}
 where the last inequality is supported by $t \geq T = \frac{4}{\mu\sqrt{s}}$. By matching the corresponding terms in~\eqref{eqn: con-lyapunov-estimate} and~\eqref{eqn: con-total-3}, we can estimate the time derivative as
    \[
     \frac{d\E}{d t} \leq - \min\left\{\frac{3\mu \sqrt{s}}{8}, \; \frac{\mu \sqrt{s}}{2}, \;\frac{2}{3\sqrt{s}},\; \frac{1}{4\sqrt{s}} \right\} \E \leq -\frac{\mu\sqrt{s}}{4}\E,
    \]
 with the latter inequality following because $ \mu \leq L$ and $s \in (0, 1/L)$.  Additionally, owing to the condition $f \in \mathcal{S}_{\mu,L}^2$, it holds that: 
 \[
  \norm{\nabla f(X)}^2 \leq 2L\brak{f(X) - f(x^\star)}.
 \]
 Hence, the proof is complete with some elementary calculations.

%    Here we use the equality $\sqrt{s}t\brak{t + \frac{\sqrt{s}}{2}} = \brak{\frac{\sqrt{s}}{2}t^2 + \frac{s}{2}t} + \frac{\sqrt{s}}{2}t^2$.
%
%
%    We need let the proportion of each corresponding coefficients larger that $\frac{\mu\sqrt{s}}{4}$. To do this, after verification, we only need to let the proportion of coefficient of $f(X) - f(x^\star)$ satisfies the condition, i.e.,
%    \begin{equation*}
%        \dfrac{\mu\sqrt{s}t^2 + (\mu s - 2)t - \sqrt{s}}{2t(t+\sqrt{s})} \ge \dfrac{\mu\sqrt{s}}{4}.
%    \end{equation*}
%    This is equivalent to 
%    \begin{equation*}
%        \frac{\mu\sqrt{s}}{2}t^2 + \brak{\frac{\mu s}{2} - 2}t - \sqrt{s} \ge 0.
%    \end{equation*}
%    This is sufficient when $t > \frac{4}{\mu\sqrt{s}}$, that is, set $T = \frac{4}{\mu\sqrt{s}}$, then we have 
%    \begin{equation*}
%        \frac{\dd}{\dd t}\E(t) \le -\frac{\mu\sqrt{s}}{4}\E(t).
%    \end{equation*}
%    The Gronwall's inequality yields the result.
\end{proof}

\section{Linear convergence via a novel discrete Lyapunov function}
\label{sec: linear-convergence}

In this section, we develop a novel discrete Lyapunov function aimed at deducing the linear convergence properties of the~\texttt{NAG} method when implemented on $\mu$-strongly convex functions. Our approach is fundamentally anchored in the theoretical framework of potential and mixed energy, as thoroughly detailed in~\cite[(5.14)]{li2022proximal} for convex scenarios. The primary distinction of our method lies in the adoption of an iteration-varying coefficient in the construction of kinetic energy, which ensures that the respective terms remain proportional. This critical alteration marks a significant departure from the techniques described in~\cite[(5.3)]{li2022linear}.

\subsection{Smooth optimization with~\texttt{NAG}}
\label{subsec: nesterov-1983}

% Taking the velocity iteration sequence as $v_{k} =(x_{k} - x_{k-1})/\sqrt{s}$, we transform~\texttt{NAG} into the implicit-velocity scheme (phase-space representation) as

In the arena of smooth optimization, we leverage the velocity iteration sequence defined by $v_{k} =(x_{k} - x_{k-1})/\sqrt{s}$, which allows us to recast~\texttt{NAG} into the implicit-velocity scheme characterized by a phase-space representation: 
\begin{equation}\label{eqn: im-vel}
   \left\{ \begin{aligned}
            & x_{k+1} - x_k = \sqrt{s}v_{k+1},\\
            & v_{k+1} - v_k = -\dfrac{r+1}{k+r}\cdot v_k - \sqrt{s}\nabla f(y_k),
           \end{aligned} \right.
\end{equation}
where the sequence $\{y_{k}\}_{k=0}^{\infty}$ complies with the relation:
\begin{equation*}
    y_k = x_k + \dfrac{k-1}{k+r}\cdot \sqrt{s}v_k.
\end{equation*}
It is evident that the second iteration of~\eqref{eqn: im-vel} can be neatly reformulated as 
\begin{equation}\label{eqn: reformulated}
    (k+r)v_{k+1} - (k-1)v_k = -(k+r)\sqrt{s}\nabla f(y_k).
\end{equation}
Building upon this framework, we now proceed to elucidate how the Lyapunov function is systematically constructed by employing the principled approach as outlined in~\cite{chen2022revisiting}. 
\begin{itemize}
\item[(\textbf{I})] Drawing upon the discrete Lyapunov function established in~\cite[(5.14)]{li2022proximal}, we here choose to implement the analogous mixed energy for its computational efficacy as referenced in~\cite[Section 4.2]{chen2022gradient}. The mixed energy is delineated as:
\begin{equation}
\label{eqn: mix-discrete}
\E_{\textbf{mix}}(k) = \frac12 \| \sqrt{s}(k-1)v_k + r(x_{k} - x^\star)  \|^2.
\end{equation}
By employing the initial step of the phase-space representation~\eqref{eqn: im-vel} and the reformulated expression provided in~\eqref{eqn: reformulated}, we establish the subsequent equality as
\begin{align}
        & \sqrt{s}kv_{k+1} + r(x_{k+1} - x^\star) -  \sqrt{s}(k-1)v_k + r(x_{k} - x^\star) \nonumber \\ 
  =    & \sqrt{s}kv_{k+1} + r(x_{k+1} - x_{k}) -  \sqrt{s}(k-1)v_k                                   \nonumber \\ 
  =    &  \sqrt{s}(k+r)v_{k+1}  -  \sqrt{s}(k-1)v_k                                                         \nonumber \\
  =    & -(k+r)s\nabla f(y_k), \label{eqn: mix-revise}
\end{align}
where the second equality is derived from the initial step of the phase-space representation~\eqref{eqn: im-vel} and the last equality follows from equation~\eqref{eqn: reformulated}. We are now in a position to determine the iterative difference of mixed energy, which is calculated as follows:
\begin{align}
     \E_{\textbf{mix}}(k+1) & - \E_{\textbf{mix}}(k) \nonumber\\ 
    & =  \big\langle\sqrt{s}kv_{k+1} + r(x_{k+1} - x^\star), -s(k+r)\nabla f(y_k)\big\rangle \nonumber \\
   & \mathrel{\phantom{=}} - \dfrac{s^2(k+r)^2}{2}\norm{\nabla f(y_k)}^2 \nonumber \\
%=  & \inner{\sqrt{s}kv_{k+1} + r(y_k - x^\star) - srG_s(y_k), -s(k+r)G_s(y_k)}- \dfrac{s^2(k+r)^2}{2}\norm{G_s(y_k)}^2\\
    & =   - \underbrace{s^{\frac{3}{2}}k(k+r)\big\langle \nabla f(y_k), v_{k+1}\big\rangle}_{\mathbf{I}} - sr(k+r)\big\langle\nabla f(y_k), y_k - x^{\star}\big\rangle \nonumber \\
    & \mathrel{\phantom{=}} - \dfrac{s^2}{2}(k-r)(k+r)\norm{\nabla f(y_k)}^2, \label{eqn: iter-mix-discrete}
\end{align}
where the ultimate step incorporates the gradient step inherent in the \texttt{NAG} method. 

\item[(\textbf{II})]  In alignment with the principled approach to constructing the Lyapunov function as outlined in~\cite{chen2022revisiting}, we commence our analysis by considering the kinetic energy. Departing from~\cite[(4.10)]{chen2022revisiting}, we formulate the kinetic energy with an iteration-varying coefficient, which is articulated as follows:
\begin{equation}
\label{eqn: kinetic-discrete}
\E_{\textbf{kin}}(k) = \frac{s\tau(k)}{2}\| v_k \|^2.
\end{equation}
By employing the expression from~\eqref{eqn: reformulated}, the iterative difference of kinetic energy can be calculated as follows:
\begin{align}
\E_{\textbf{kin}}(k+1) &- \E_{\textbf{kin}}(k)  \nonumber \\ %  =  \frac{s\tau(k+1)}{2}\norm{v_{k+1}}^2 - \frac{s\tau(k)}{2}\norm{v_k}^2  \nonumber \\
                                              & =  \frac{s\tau(k+1)}{2}\norm{v_{k+1}}^2 - \frac{s\tau(k)}{2}\left( \frac{k + r}{k - 1} \right)^2\norm{v_{k+1} + \sqrt{s}\nabla f(y_k)}^2 \nonumber \\
                                              & =  \frac{s}{2} \bigg[ \tau(k+1) - \tau(k)\left( \frac{k + r}{k - 1} \right)^2 \bigg]\norm{v_{k+1}}^2 \nonumber \\
                                              &\mathrel{\phantom{=}}  - s^{\frac32}\tau(k)\left( \frac{k + r}{k - 1} \right)^2 \big\langle\nabla f(y_k), v_{k+1}\big\rangle \nonumber \\
                                              & \mathrel{\phantom{=}} - \frac{s^2  \tau(k)}{2}\left( \frac{k + r}{k - 1} \right)^2\norm{\nabla f(y_k)}^2  . \label{eqn: iter-kin-discrete1}
    \end{align}

To ease computation as indicated by~\eqref{eqn: iter-kin-discrete1}, an intuitive choice for $\tau(k)$ would be $\tau(k) = (k-1)^2$. This selection simplifies the iterative difference~\eqref{eqn: iter-kin-discrete1}, which can be succinctly expressed as
\begin{align}
\E_{\textbf{kin}}(k+1) - \E_{\textbf{kin}}(k) = & - \frac{sr(2k + r)}{2} \|v_{k+1}\|^2  - \frac{s^2  (k+r)^2}{2} \norm{\nabla f(y_k)}^2 \nonumber \\  
                                                                       & - \underbrace{s^{\frac32}( k + r )^2 \big\langle\nabla f(y_k), v_{k+1}\big\rangle}_{\mathbf{II}}.  \label{eqn: iter-kin-discrete}
\end{align}

%Obviously, a simple way to choose the coefficient is $\tau(k) = (k-1)^2$.

\item[(\textbf{III})] We now turn our attention to the potential energy.  It is appropriate for the potential energy to also feature an iteration-varying coefficient, hereafter referred to as $\gamma(k)$. Consequently,  we define the potential energy as: 
\begin{equation}
\label{eqn: pot-discrete}
\E_{\textbf{pot}}(k) = s\gamma(k) \left( f(x_k) - f(x^\star) \right).
\end{equation}
To determine the iterative difference in potential energy, we must refer back to the inequality for $f \in \mathcal{S}_{\mu,L}^1$, as demonstrated in~\cite{nesterov1998introductory}:  
\begin{align}
f(y - s\nabla f(y)) & - f(x) \nonumber \\ &\leq \big\langle \nabla f(y), y - x \big\rangle - \frac{\mu}{2}\|y - x\|^2 - \left( s - \frac{Ls^2}{2} \right) \| \nabla f(y) \|^2. \label{eqn: prox-smooth}
\end{align}
Applying~\eqref{eqn: prox-smooth} with $x_{k+1}$ and $x_{k}$, we derive:
\begin{align}
 f(x_{k+1}) & -  f(x_k) \nonumber \\
                 & \leq  \big\langle\nabla f(y_k), y_k - x_k \big\rangle  - \frac{\mu}   {2}\norm{y_k - x_k}^2 - \brak{s - \frac{s^2 L}{2}}\norm{\nabla f(y_k)}^2 \nonumber \\
                 & =  \big\langle\nabla f(y_k), \sqrt{s}v_{k+1} +   s\nabla f(y_k)\big\rangle - \frac{\mu}{2}\norm {\sqrt{s}v_{k+1} + s\nabla f(y_k)}^2 \nonumber  \\
                 & \mathrel{\phantom{=}} - \brak{s - \frac{s^2 L}{2}}   \norm{\nabla f(y_k)}^2  \nonumber \\
%            = & \sqrt{s}\inner{G_s(y_k), v_{k+1}} +  \frac{s}{2}\norm{G_s(y_k)}^2 - \dfrac{s  (1-sL)}{2}\norm{G_s(y_k)}^2 - \frac{\mu s}{2}\norm{v_{k+1}}^2 \nonumber \\
%            &- \frac{\mu s^2}{2}\norm{G_s(y_k)}^2 -    \mu s^{\frac{3}{2}}\inner{G_s(y_k), v_{k   +1}} \nonumber\\
                  & =  \sqrt{s}(1-\mu s)\big\langle\nabla f(y_k), v_ {k+1}\big\rangle + \frac{s^2(L - \mu)}{2}\norm{\nabla f   (y_k)}^2 -  \frac{\mu s}{2}\norm{v_{k+1}}^2 \label{eqn: key1}
\end{align}
where the first equality is justified by the initial iteration in~\eqref{eqn: im-vel} in conjunction with the gradient step of~\texttt{NAG}.  Subsequently, for the potential energy~\eqref{eqn: pot-discrete}, the iterative difference can be estimated as follows:
\begin{align}
& \E_{\textbf{pot}}(k+1)   - \E_{\textbf{pot}}(k) \nonumber \\
& =  s\gamma(k+1) \left( f(x_{k+1}) - f(x^\star) \right) -  s\gamma(k) \left( f(x_k) - f(x^\star) \right)\nonumber  \\
& =  s\gamma(k) \left( f(x_{k+1}) - f(x_{k}) \right) + s (\gamma(k+1) - \gamma(k)) \left( f(x_{k+1}) - f(x^\star) \right) \nonumber  \\
& \leq  s\gamma(k) \left(\underbrace{\sqrt{s}(1-\mu s)\big\langle\nabla f(y_k), v_ {k+1}\big\rangle}_{\mathbf{III}} + \frac{s^2(L - \mu)}{2}\norm{\nabla f   (y_k)}^2 -  \frac{\mu s}{2}\norm{v_{k+1}}^2 \right)\nonumber  \\
& \mathrel{\phantom{=====}} + \; s (\gamma(k+1) - \gamma(k)) \left(  f(x_{k+1}) - f(x^\star) \right). \label{eqn: iter-pot-discrete}
\end{align}
Coupled with the earlier derivations,~\eqref{eqn: iter-mix-discrete} and~\eqref{eqn: iter-kin-discrete}, it becomes imperative to fine-tune the coefficient $\gamma(k)$ within the iterative difference~\eqref{eqn: iter-pot-discrete} to cancel out the term $\big\langle\nabla f(y_k), v_ {k+1}\big\rangle$. Explicitly, this necessitates the fulfillment of the relation:
\[
s\gamma(k)\cdot\mathbf{III} - \mathbf{I} - \mathbf{II} = 0.
\]
Therefore, the iteration-varying coefficient $\gamma(k)$ should be determined as follows:
\[
\gamma(k) = \frac{(k + r)(2k + r)}{1 - \mu s}. 
\]
\end{itemize}

Integrating the mixed energy~\eqref{eqn: mix-discrete}, the kinetic energy~\eqref{eqn: kinetic-discrete} and the potential energy~\eqref{eqn: pot-discrete} --- each bolstered by iteration-varying coefficients --- we construct the noval discrete Lyapunov function as depicted below: 
\begin{align}
\E(k) =  &\dfrac{s(k+r)(2k + r)}{1-\mu s}\brak{f(x_k) - f(x^\star)} \nonumber \\ & + \frac{s(k-1)^2}{2}\norm{v_k}^2 + \frac{1}{2}\norm{\sqrt{s}(k-1)v_k + r(x_k - x^\star)}^2. \label{eqn: lyapunov-dis}
\end{align}
Within this novel Lyapunov function~\eqref{eqn: lyapunov-dis} at our disposal, we methodically infer the convergence rates for both the function value and the squared gradient norm, as stated in the forthcoming theorem. 

%With~\eqref{eqn: lyapunov-dis}, we conclude the section with the following theorem and its proof.
\begin{theorem}
\label{thm: prox-smooth}
 Let $f \in \mathcal{S}_{\mu,L}^{1}$. Given any step size $0 < s < 1/L$, there exists a positive integer $K := K(L, \mu, s, r) $ such that the iterative sequence $\{(x_{k}, y_{k})\}_{k=0}^{\infty}$ generated by~\texttt{NAG} with any initial $x_{0} = y_{0}\in \mathbb{R}^d$ satisfies 
\begin{equation}
\label{eqn: error-convergence}      
\left\{ \begin{aligned}
             f(x_k) - f(x^\star) & \leq   \dfrac{\E(K)}{s(k+r)(2k+r)\left[1 + (1-Ls) \cdot \dfrac{\mu s}{4}\right]^{k-K}} \\
         \|\nabla f(y_{k})\|^2 & \leq \dfrac{4\E(K)}{s^2(1 - Ls)(k+r)(2k+r)\left[1 + (1-Ls) \cdot \dfrac{\mu s}{4}\right]^{k-K}}
         \end{aligned} \right.
\end{equation}
for any $k \geq K$. 
\end{theorem}

\begin{remark}
Considering the inverse proportionality between the step size and the Lipshitz constant, denoted as $s \sim 1/L$,  we articulate the convergence rates specified in~\eqref{eqn: error-convergence} in terms of the dimensionless parameter $\alpha = sL $, which falls within the open interval $(0,1)$. Accordingly, the established bounds may be recast as 
\begin{equation}
\label{eqn: error-convergence-1}      
\left\{ \begin{aligned}
        f(x_k) - f(x^\star) &\leq \frac{\E(K)}{(k+r)(2k+r)\brak{1 + \frac{\alpha(1-\alpha)}{4} \cdot \frac{\mu }{L}}^{k-K}} \\
        \|\nabla f(y_{k})\|^2 &\leq \dfrac{4L^2\E(K)}{\alpha^2(1 - \alpha)(k+r)(2k+r)\left[1 + \frac{\alpha(1 - \alpha)}{4} \cdot \dfrac{\mu}{L}\right]^{k-K}}
          \end{aligned} \right.
\end{equation}
for any $k \geq K$. The bounds delineated in~\eqref{eqn: error-convergence-1} indicate that both the function value and the square of the gradient norm adhere to a linear convergence pattern, characterized by a rate akin to $\left( 1 + c\cdot\frac{\mu}{L}\right)^{-k}$ where the constant $c$ resides in the range $(0,1)$. Albeit the constant $c$ can be optimized to $1$, the ensuing convergence rate of $\left(1 +\frac{\mu}{L}\right)^{-k}$ still exhibits a significant gap when compared to the optimal rates of $\left(1 - \sqrt{\frac{\mu}{L}}\right)^{k}$ or $\left( 1 + \sqrt{\frac{\mu}{L}} \right)^{-k}$, as elucidated in~\cite{nesterov1998introductory}.  This discrepancy accentuates the sustained pursuit in the realm of computational optimization for devising accelerated methods that circumvent any prerequisites concerning the modulus of strong convexity.
\end{remark}

%As we know, the scale of the step size is $s \sim 1/L$. Thus, we can rewrite the convergence rate~\eqref{eqn: error-convergence} in accordance with the big-O form by the nondimensional parameter $\alpha = sL \in (0,1)$.
%, we can obtain the iterative difference of the total Lyapunov function as
%  the calculation for these three energy functions above, we have 
%       Next recall that applying Cauchy-Schwarz inequality, the Lyapunov function satisfies 
%\[
%\|\nabla f(y_{k})\|^2 \leq \dfrac{4\E(K)}{s^2(1 - Ls)(k+r)(2k+r)\left[1 + (1-Ls) \cdot \dfrac{\mu s}{4}\right]^{k-K}}
%\]

\begin{proof}[Proof of~\Cref{thm: prox-smooth}]
Given the discrete Lyapunov function defined in~\eqref{eqn: lyapunov-dis}, we calculate its iterative difference by summing the three iterative differences laid out in~\eqref{eqn: iter-mix-discrete},~\eqref{eqn: iter-kin-discrete} and~\eqref{eqn: iter-pot-discrete}.  This amalgamation yields:
    \begin{align}
                       \E(k+1) & - \E(k) \nonumber \\
                      & = - \underbrace{sr(k+r)\big\langle \nabla f(y_k), y_k - x^\star\big\rangle}_{\mathbf{IV}}  - \frac{sr(2k+r)}{2}\norm{v_{k+1}}^2 \nonumber \\
                      & \mathrel{\phantom{=}} - s^2k(k+r)\norm{\nabla f(y_k)}^2 \nonumber \\
                      & \mathrel{\phantom{=}} + \dfrac{s(k+r)(2k + r)}{1-\mu s}\brak{ \frac{s^2(L - \mu)}{2}\norm{\nabla f   (y_k)}^2 -  \frac{\mu s}{2}\norm{v_{k+1}}^2 } \nonumber  \\
                      & \mathrel{\phantom{=}} + \frac{s(4k + 3r + 2)}{1-\mu s}\brak{ f(x_{k+1}) - f(x^{\star}) } \label{eqn: total-iter-diff}
    \end{align}
To ensure proportional alignment of terms, we exhibit $\|x_{k+1} - x^{\star} \|^2$ and the other terms in~\eqref{eqn: total-iter-diff}, thereby replacing Term $\mathbf{IV}$. After inserting $x_{k+1}$ and $x^{\star}$ into~\eqref{eqn: prox-smooth}, we deduce:
    \begin{align} 
        f(x_{k+1}) - f(x^\star) &\leq \big\langle\nabla f(y_k), y_k - x^\star\big\rangle - \dfrac{s}{2}\norm{\nabla f(y_k)}^2 - \dfrac{\mu}{2}\norm{y_k - x^\star}^2 \nonumber \\
                                          &\leq \big\langle\nabla f(y_k), y_k - x^\star\big\rangle - \dfrac{s}{2}\norm{\nabla f(y_k)}^2 - \dfrac{\mu}{2}\norm{x_{k+1} - x^\star}^2, \label{eqn: key2}
    \end{align}
    where the penultimate step is intrinsically connected to the gradient method, such that:
    \begin{align*}
        \norm{x_{k+1} - x^\star}^2 = & \norm{y_k - x^\star - s \nabla f(y_k)}^2 \\
        = & \norm{y_k - x^\star}^2 - 2s\big\langle\nabla f(y_k), y_k - x^\star\big\rangle + s^2\norm{\nabla f(y_k)}^2  \leq \norm{y_k - x^\star}^2,
    \end{align*}
    where the last inequality holds for any $s < 1/L$ due to the first inequality of~\eqref{eqn: key2}. By substituting inequality~\eqref{eqn: key2} into the iterative difference~\eqref{eqn: total-iter-diff}, we attain:  
    \begin{align}
        \E(k+1) - \E(k) \le & s\brak{\frac{4k + 3r + 2}{1-\mu s} - r(k+r)}\brak{f(x_{k+1}) - f(x^\star)} \nonumber \\
                                     & - \frac{s}{2}\brak{\frac{\mu s}{1-\mu s}(k+r)(2k + r) + r(2k + r)}\norm{v_{k+1}}^2 \nonumber\\
                                     & - \frac{\mu s }{2} \cdot r(k+r)\norm{x_{k+1} - x^\star}^2 \nonumber \\
                                     &  - \frac{s^2}{2} \cdot \frac{1-Ls}{1-\mu s}  \cdot (k+r)(2k + r) \norm{\nabla f(y_k)}^2. \label{eqn: total-iter-diff2}
    \end{align}
   
    To proportionally adjust the corresponding terms within the Lyapunov function $\mathcal{E}(k+1)$, we employ the Cauchy-Schwarz inequality to estimate as follows:
   \begin{align}
      \E(k+1)  &\leq \frac{s(k+r+1)(2k+r+2)}{1-\mu s}(f(x_{k+1}) - f(x^\star)) \nonumber \\
                   & \mathrel{\phantom{\leq}}  + \frac{3}{2}sk^2\norm{v_{k+1}}^2 + r^2\norm{x_{k+1} - x^\star}^2 \nonumber\\
                    &= \frac{s\big[(k+r)(2k+r) + (4k+3r+2)\big]}{1-\mu s}(f(x_{k+1}) - f(x^\star)) \nonumber\\
                    &\mathrel{\phantom{=}}+ \frac{3}{2}sk^2\norm{v_{k+1}}^2 + r^2\norm{x_{k+1} - x^\star}^2. \label{eqn: lyapunov-estimate}
   \end{align}
For any $f \in \mathcal{S}_{\mu,L}^{1}$, the following inequality is satisfied:
\begin{equation}
\label{eqn: gradient-strong-general}
\|\nabla f(y)\|^2 \geq 2\mu\left( f(y - s \nabla f(y)) - f(x^{\star}) \right)
\end{equation}
 for any $y \in \mathbb{R}^{d}$.  Applying this to $y_{k}$ of~\texttt{NAG} within~\eqref{eqn: gradient-strong-general}, we arrive at: 
     \begin{equation}
    \label{eqn: gradient-strong}
     \|\nabla f(y_k)\|^2 \geq 2\mu\left( f(x_{k+1}) - f(x^{\star}) \right).
    \end{equation}
Incorporating inequality~\eqref{eqn: gradient-strong} into the iterative difference~\eqref{eqn: total-iter-diff2}, we deduce: 
    \begin{align}
        \E(k+1) & - \E(k) \nonumber \\
        \leq & -\frac{s}{1-\mu s}\brak{\frac{\mu s(1 - Ls)}{2}(k+r)(2k+r) - (4k+3r+2) + r(1-\mu s)(k+r)} \nonumber \\
                & \cdot \brak{f(x_{k+1}) - f(x^\star)} - \frac{\mu s}{2}r(k+r)\norm{x_{k+1} - x^\star}^2  \nonumber \\
        & - \frac{s}{2}\brak{\frac{\mu s}{1-\mu s}(k+r)(2k + r) + r(2k + r)}\norm{v_{k+1}}^2 \nonumber \\
        & - \frac{s^2}{4}\frac{1-Ls}{1-\mu s}(k+r)(2k + r)\norm{\nabla f(y_k)}^2 \nonumber \\
        \leq & -\frac{s}{1-\mu s}\brak{\frac{\mu s(1 - Ls)}{2}(k+r)(2k+r) - (4k+3r+2) + r(1-\mu s)(k+r)} \nonumber \\
               & \cdot \brak{f(x_{k+1}) - f(x^\star)}  - \frac{\mu sr^2}{2}\norm{x_{k+1} - x^\star}^2 \nonumber \\
        & - \mu s^2k^2\norm{v_{k+1}}^2 - \frac{s^2}{4}(1-Ls)(k+r)(2k + r)\norm{\nabla f(y_k)}^2 \label{eqn: total-diff3}
    \end{align}
To align the terms between the right-hand side of~\eqref{eqn: total-diff3} and~\eqref{eqn: lyapunov-estimate}, we establish the following inequality as
   \begin{multline*}
    (1-Ls) \cdot \frac{\mu s}{4}(k+r)(2k+r) - (4k+3r+2) + r(1-\mu s)(k+r) \\ \geq  \frac{\mu s(1-Ls) }{4} \cdot (4k+3r+2).
   \end{multline*}
   It is evident that there exists a positive constant $K = K(L, \mu, s, \alpha )$ for which the above inequality holds.  Thus, we conclude:
   \begin{align*}
                   & \E(k+1) -\E(k) \\
                   & \le - \min\left\{\frac{1-Ls}{4},\; \frac{2}{3},\; \frac{1}{2} \right\} \mu s \cdot \E(k+1) - \frac{s^2(1-Ls)(k+r)(2k+r)}{4}\norm{\nabla f(y_k)}^2 \\
                   & = \frac{1-Ls}{4} \cdot \mu s \cdot \E(k+1) - \frac{s^2(1-Ls)(k+r)(2k+r)}{4}\norm{\nabla f(y_k)}^2.
   \end{align*}
   The proof is complete with some elementary operations.
\end{proof}

\subsection{Composite optimization via~\texttt{FISTA}}
\label{subsec: fista}

In this section, we extend the convergence rates of~\texttt{NAG} as established in~\Cref{thm: prox-smooth} to include its proximal variant,~\texttt{FISTA}. As specified in~\Cref{defn: proximal-subgradient},~\texttt{FISTA} utilizes the $s$-proximal operator~\eqref{eqn: proximal-operator} and is defined by the following iterative scheme, starting from any initial point $y_0 = x_0 \in \mathbb{R}^d$:
\[
\left\{ \begin{aligned}
        & x_{k} = P_s\left(y_{k-1} - s\nabla f(y_{k-1})\right), \\
        & y_{k} = x_{k} + \frac{k-1}{k+r}(x_k - x_{k-1}),
        \end{aligned} \right.
\]
where $s > 0$ denotes the step size. We can analogously formulate~\texttt{FISTA} in a manner akin to~\texttt{NAG} using the $s$-proximal subgradient operator~\eqref{eqn: subgradient-operator}, which yields
\begin{equation}
\label{eqn: fista-nag}
\left\{ \begin{aligned}
        & x_{k} = y_{k-1} - sG_s(y_{k-1}), \\
        & y_{k} = x_{k} + \frac{k-1}{k+r}(x_k - x_{k-1}),
        \end{aligned} \right.
\end{equation}
where the proximal operator $G_s(\cdot)$ replace the gradient operator  $\nabla f(\cdot)$ in~\texttt{NAG}.  Furthermore, by designing the velocity iteration sequence as $v_{k} =(x_{k} - x_{k-1})/\sqrt{s}$, we reformulate the~\texttt{FISTA} updates in a \texttt{NAG}-esque fashion~\eqref{eqn: fista-nag} into a implicit-velocity scheme (phase-space representation), expressed as
\begin{equation}
\label{eqn: im-vel-fista}
   \left\{ \begin{aligned}
            & x_{k+1} - x_k = \sqrt{s}v_{k+1},\\
            & v_{k+1} - v_k = -\dfrac{r+1}{k+r}\cdot v_k - \sqrt{s} G_s(y_k).
           \end{aligned} \right.
\end{equation}

To derive the convergence rates, we still need to establish a discrete Lyapunov function. Here, we generalize the one previously utilized for smooth functions~\eqref{eqn: lyapunov-dis}, by substituting the smooth function $f$ with the composite function $\Phi = f+ g$, resulting in
\begin{align}
\E(k) =  &\dfrac{s(k+r)(2k + r)}{1-\mu s}\brak{\Phi(x_k) - \Phi(x^\star)} \nonumber \\ & + \frac{s(k-1)^2}{2}\norm{v_k}^2 + \frac{1}{2}\norm{\sqrt{s}(k-1)v_k + r(x_k - x^\star)}^2. \label{eqn: lyapunov-dis-fista}
\end{align}
As outlined in~\Cref{subsec: nesterov-1983}, deriving the desired convergence rates for the smooth case predominantly hinges on two key inequalities,~\eqref{eqn: prox-smooth} and~\eqref{eqn: gradient-strong-general}.  If we can successfully adjust these two key inequalities to account for the proximal setting, we shall be able to transpose the results of~\Cref{thm: prox-smooth} into this broader context. For the fundamental inequality~\eqref{eqn: prox-smooth},  its proximal analogue is articulated as 
\begin{align}
\Phi(y - sG_s(y)) & - \Phi(x) \nonumber \\ &\leq \big\langle G_s(y), y - x \big\rangle - \frac{\mu}{2}\|y - x\|^2 - \left( s - \frac{Ls^2}{2} \right)  \| G_s(y) \|^2, \label{eqn: prox-composite}
\end{align}
where the proof can be found in~\cite[Lemma 3.2]{li2022linear}. The primary challenge lies in extending the $\mu$-strongly convex inequality~\eqref{eqn: gradient-strong-general} to the composite function $\Phi = f+ g$, for which we establish the subsequent lemma.

\begin{lemma} \label{lem: proximal-key-inequality}
    Let $f \in \mathcal{S}_{\mu,L}^{1}$ and $g \in \mathcal{F}^{0}$. It then holds that the $s$-proximal subgradient, as defined in~\eqref{eqn: subgradient-operator}, satisfies the following inequality: 
    \begin{equation}
    \label{eqn: proximal-key-inequality}   
        \norm{G_s(y)}^2 \ge 2\mu\brak{\Phi(y - sG_s(y)) - \Phi(x^\star)},
    \end{equation}
    for any $y \in \mathbb{R}^d$.
\end{lemma}

\begin{proof}[Proof of~\Cref{lem: proximal-key-inequality}]
    Under the condition that the step size adheres to $0 < s \leq 1/L $, we invoke the fundamental proximal inequality~\eqref{eqn: prox-composite}, which simplifies to:
    \begin{equation}
    \label{eqn: fund-prox-strong} 
        \Phi(y - sG_s(y)) - \Phi(x) \le \big\langle G_s(y), y - x\big\rangle - \dfrac{s}{2}\norm{G_s(y)}^2 - \frac{\mu}{2}\norm{y - x}^2
    \end{equation}
    for any $x,y \in \mathbb{R}^d$. By reorganizing the terms of~\eqref{eqn: fund-prox-strong},  we arrive at the expression 
    \begin{equation}\label{eqn: lem-prox-1}
       \Phi(x) \geq  \Phi(y - sG_s(y)) + \big\langle G_s(y), x - y\big\rangle + \frac{s}{2}\norm{G_s(y)}^2 + \frac{\mu}{2}\norm{y -x}^2 
    \end{equation}
for any $x,y \in \mathbb{R}^d$. For succinctness, we denote the right-hand side of~\eqref{eqn: lem-prox-1} as
    \begin{equation}\label{eqn: lem-prox-2}
        h(x) = \Phi(y - sG_s(y)) + \big\langle G_s(y), x - y\big\rangle + \frac{s}{2}\norm{G_s(y)}^2 + \frac{\mu}{2}\norm{y -x}^2,
    \end{equation}
for any fixed $y \in \mathbb{R}^d$. Consequently, inequality~\eqref{eqn: lem-prox-1} takes the form $\Phi(x) \geq h(x)$. Considering that $x^{\star}$ is the unique minimizer of  the composite function $\Phi$,  substituting $x$ with $x^{\star}$ yields:
    \begin{equation}\label{eqn: lem-prox-3}
       \Phi(x^{\star}) \geq  h(x^{\star}) .
    \end{equation}
Upon evaluating expression~\eqref{eqn: lem-prox-2}, we identify that $h(x)$ embodies a quadratic function whose Hessian is positive definte.   This is depicted as:
\begin{equation}\label{eqn: lem-prox-4}
h(x) = \frac{\mu}{2} \left\|x - y + \frac{1}{\mu} G_s(y) \right\|^2 + \Phi(y-sG_s(y)) - \brak{\frac{1}{2\mu} - \frac{s}{2}}\norm{G_s(y)}^2.
\end{equation}
Incorporating~\eqref{eqn: lem-prox-4} into~\eqref{eqn: lem-prox-3} results in
\begin{equation}\label{eqn: lem-prox-5}
\Phi(x^{\star}) \geq \Phi(y-sG_s(y)) - \brak{\frac{1}{2\mu} - \frac{s}{2}}\norm{G_s(y)}^2 \geq \Phi(y-sG_s(y))  - \frac{1}{2\mu} \norm{G_s(y)}^2 . 
\end{equation}
By rearranging~\eqref{eqn: lem-prox-5}, we complete the proof. 

%
%    Thus we have 
%    \begin{equation}\label{eqn: 2}
%        \Phi(y - sG_s(y)) + \inner{G_s(y), x^\star - y} + \frac{s}{2}\norm{G_s(y)}^2 + \frac{\mu}{2}\norm{y - x^\star}^2 \le \Phi(x^\star).
%    \end{equation}
%    Now we compute the minimum of the left hand side of (\ref{eqn: 1}). Let 
%
%    Finding the minimum point is equivalent to the derivative of left hand side is $0$, i.e., 
%    \begin{equation*}
%        \nabla h(x) = 0.
%    \end{equation*}
%    This is solved as 
%    \begin{equation*}
%        x - y = -\frac{1}{\mu} G_s(y)
%    \end{equation*}
%    Substituting this into (\ref{eqn: 3}), we have 
%    \begin{equation}
%        h(x) \ge \Phi(y-sG_s(y)) - \brak{\frac{1}{2\mu} - \frac{s}{2}}\norm{G_s(y)}^2, \quad \forall x\in \mr^d.
%    \end{equation}
%    Hence 
%    \begin{equation*}
%        h(x^\star)\ge \Phi(y-sG_s(y)) - \brak{\frac{1}{2\mu} - \frac{s}{2}}\norm{G_s(y)}^2
%    \end{equation*}
%    Combining this inequality with (\ref{eqn: 2}) and rearranging yields the result.
\end{proof}

With the proximal generalization of the $\mu$-strongly convex inequality, as demonstrated in~\Cref{lem: proximal-key-inequality}, we now present the subsequent theorem,  which characterizes the convergence rates of~\texttt{FISTA}. 
\begin{theorem}\label{thm: prox-fista}
    Let $f \in \mathcal{S}_{\mu,L}^{1}$ and $g \in \mathcal{F}^{0}$. For any step size $0 < s < 1/L$, there exists some positive integer $K := K(L,\mu, s, r) $, such that the iterative sequence $\{(x_{k}, y_{k})\}_{k=0}^{\infty}$ generated by~\texttt{FISTA}, with any initial $x_{0} = y_0 \in \mathbb{R}^d$, satisfies the following inequalities as
\begin{equation}
\label{eqn: fista-convergence}      
\left\{ \begin{aligned}
            \Phi(x_k) - \Phi(x^\star) &\leq   \dfrac{\E(K)}{s(k+r)(2k+r)\left[1 + (1-Ls) \cdot \dfrac{\mu s}{4}\right]^{k-K}} \\
         \|G_s(y_{k})\|^2 &\leq \dfrac{4\E(K)}{s^2(1 - Ls)(k+r)(2k+r)\left[1 + (1-Ls) \cdot \dfrac{\mu s}{4}\right]^{k-K}}
         \end{aligned} \right.
\end{equation}
%\[
%\|\nabla f(y_{k})\|^2 \leq \dfrac{4\E(K)}{s^2(1 - Ls)(k+r)(2k+r)\left[1 + (1-Ls) \cdot \dfrac{\mu s}{4}\right]^{k-K}}
%\]
for any $k \geq K$. 
\end{theorem}

%
%As we have claimed in~\Cref{sec: intro}, from~\Cref{sec: high-resolution}, the proof of convergence rate of high-resolution equation provides us enough insight on the problem. This leads us to implement a resembled discrete version of Lyapunov function to get the linear convergence rate for FISTA under $\mu$-strongly convexity. But note that, the center of the proof of continuous case is the Polyak-Lojasiewicz (LP) inequality, with whom we get able to control the growth rate of the coefficient of $f(X) - f(x^\star)$. This guides us to find a proximal conterpart of LP inequality and give rise to the following lemma.
%
%~\cite[(3.5)]{li2022linear}, 
%\[
%\Phi(y - sG_s(y)) - \Phi(x) \leq \langle G_s(y), y - x \rangle - \frac{\mu}{2}\|y - x\|^2 - \left( s - \frac{Ls^2}{2} \right) \| G_s(y) \|^2
%\]

\section{Discussion and further works}
\label{sec: conclusion-further-works}

Upon reassessment of the proofs involving the continuous high-resolution ODE as shown in~\Cref{sec: high-resolution}, we recognize the shortcomings of the implicit-velocity phase-space representation.  Instead, we find it necessary to utilize the gradient-correction phase-space representation.   The strategy for conceiving the Lyapunov function~\eqref{eqn: continuous_Lyapunov}  hinges on the exclusion of terms similar to~$\langle \nabla f(X), \dot{X} \rangle$. Furthermore, the gradient-correction term, as indicated in~\eqref{eqn: time-derive-kin}, introduces a quadratic form of the velocity $\dot{X}$, which is closely associated with the Hessian. This crucial link is indispensable for achieving linear convergence within the gradient-correction high-resolution ODE (3.1), especially for $\mu$-strongly convex functions, where the minimum eigenvalue of the Hessian is bounded by the parameter $\mu$. Conversely, as elucidated in~\Cref{sec: linear-convergence},  the implicit-velocity phase-space representation is identified as the preferred method for the discrete algorithms like~\texttt{NAG} and~\texttt{FISTA}, as substantiated by references~\cite{chen2022revisiting,li2022linear}.  A notable discrepancy emerges between the continuous high-resolution ODE and the discrete algorithms, such as~\texttt{NAG} and~\texttt{FISTA}.  Nevertheless, the continuous high-resolution ODE  still yields insights into the linear convergence for these discrete algorithms, augmented by the concept of implicit velocity, which necessitates the strategic exclusion of terms related to $\langle\nabla f(y_k), v_ {k+1}\rangle$. The utilization of the implicit-velocity phase-space representation offers two considerable benefits: it significantly simplifies computational processes and facilitates the construction of the Lyapunov function~\eqref{eqn: lyapunov-dis}.

%This choice offers numerous advantages, primarily in significantly simplifying calculations. The construction of the Lyapunov function~\eqref{eqn: lyapunov-dis} is facilitated and straightforward when we exclude terms associated with $\langle\nabla f(y_k), v_ {k+1}\rangle$. 

% Conversely, for the continuous high-resolution ODE, the implicit-velocity phase representation is unsuitable. We are instead compelled to employ the gradient-correction one. Despite this, the philosophy of eliminating terms akin to~$\big\langle \nabla f(X), \dot{X} \big\rangle$ when formulating the Lyapunov function~\eqref{eqn: continuous_Lyapunov} remains consistent. 

%, the implicit-velocity one does not work instead but we have to adopt the gradient-correction one, although it is consistent to eliminate the terms related to  $\big\langle \nabla f(X), \dot{X} \big\rangle$ for the construction of the Lyapunov function~\eqref{eqn: continuous_Lyapunov}. In addition, we can also find from~\eqref{eqn: time-derive-kin} that the gradient-correction term brings about the quadratic form of the velocity $\dot{X}$ associated with the Hessian. Indeed, this is the key reason why we can obtain the linear convergence of the gradient-correction high-resolution ODE~\eqref{eqn: grad-correction-ode} on the $\mu$-strongly convex function, since the minimum eigenvalue of the Hessian is bounded by the parameter $\mu$.
%
% In this respect, it is not clear in the proof of the discrete~\texttt{NAG} shown in~\Cref{sec: linear-convergence}. Hence, we conclude this section with two open questions. 

\begin{tcolorbox}
\begin{itemize}
\item Is it feasible to establish consistent proofs built upon a unified phase-space representation for the discrete algorithms, such as~\texttt{NAG} and~\texttt{FISTA}, and the continuous high-resolution ODE, in the context of $\mu$-strongly convex functions? 
\item Can the exposition of linear convergence for~\texttt{NAG} and~\texttt{FISTA}, as showcased in~\Cref{sec: linear-convergence}, be enhanced to yield a profound intuitive understanding similar to that offered by the gradient-correction high-resolution ODE?
\end{itemize}
\end{tcolorbox}

\bibliographystyle{siamplain}
\bibliography{reference}

\end{document}